\pdfoutput=1

\makeatletter
\IfFileExists{./numapde-latex/numapde-packages.sty}{\providecommand*{\input@path}{}\edef\input@path{{./numapde-latex/}\input@path}}{}
\makeatother

\documentclass{numapde-preprint}

\usepackage{numapde-semantic}
\usepackage{numapde-style}
\usepackage{numapde-local}

\IfFileExists{./numapde-bibliography/numapde.bib}{\addbibresource{./numapde-bibliography/numapde.bib}}{\addbibresource{numapde.bib}}

\hypersetup{
	pdftitle={First- and Second-Order Analysis for Optimization Problems with Manifold-Valued Constraints},
	pdfauthor={Ronny Bergmann, Roland Herzog, Julián Ortiz López, Anton Schiela},
	pdfkeywords={optimization on manifolds, manifold-valued constraints, manifold with corners, first- and second-order optimality conditions, Lagrangian function}
}

\title{First- and Second-Order Analysis for Optimization Problems with Manifold-Valued Constraints\thanks{This work was supported by DFG grants SCHI~1379/3--1 as well as HE~6077/10--1 within the \href{https://spp1962.wias-berlin.de}{Priority Program SPP~1962} (\emph{Non-smooth and Complementarity-based Distributed Parameter Systems: Simulation and Hierarchical Optimization}), which is gratefully acknowledged.}}
\subtitle{}
\shorttitle{Optimization Problems with Manifold-Valued Constraints}

\author{Ronny Bergmann\thanks{Norwegian University of Science and Technology, Department of Mathematical Sciences, NO-7491 Trondheim, Norway (\email{ronny.bergmannn@ntnu.no}, \url{https://www.ntnu.edu/employees/ronny.bergmann}, \orcid{0000-0001-8342-7218}).}
\and
Roland Herzog\thanks{Interdisciplinary Center for Scientific Computing, Heidelberg University, 69120 Heidelberg, Germany (\email{roland.herzog@iwr.uni-heidelberg.de}, \url{https://scoop.iwr.uni-heidelberg.de}, \orcid{0000-0003-2164-6575}).}
\and
Julián Ortiz López\thanks{Department of Mathematics, University of Bayreuth, 95440 Bayreuth, Germany (\email{julian.ortiz-lopez@uni-bayreuth.de}, \email{anton.schiela@uni-bayreuth.de}, \url{https://num.math.uni-bayreuth.de/en/team/anton-schiela/}, \orcid{0000-0002-6959-2951}).}
\and
Anton Schiela\footnotemark[4]}
\shortauthor{R. Bergmann, R. Herzog, J. Ortiz López and A. Schiela}

\dedication{}

\IfFileExists{numapde-uncolorizeHyperrefIfChanges.sty}{\RequirePackage{numapde-uncolorizeHyperrefIfChanges}}{}

\begin{document}
\maketitle

\begin{abstract}
We consider optimization problems with manifold-valued constraints.
These generalize classical equality and inequality constraints to a setting in which both the domain and the codomain of the constraint mapping are smooth manifolds.
We model the feasible set as the preimage of a submanifold with corners of the codomain.
The latter is a subset which corresponds to a convex cone locally in suitable charts.
We study first- and second-order optimality conditions for this class of problems.
We also show the invariance of the relevant quantities with respect to local representations of the problem.\end{abstract}


\begin{AMS}
\href{https://mathscinet.ams.org/msc/msc2010.html?t=90C30}{90C30}, \href{https://mathscinet.ams.org/msc/msc2010.html?t=90C46}{90C46}, \href{https://mathscinet.ams.org/msc/msc2010.html?t=49Q99}{49Q99}, \href{https://mathscinet.ams.org/msc/msc2010.html?t=65K05}{65K05}
\end{AMS}

\section{Introduction}
\label{section:Introduction}

The presence of constraints renders optimization problems not only more interesting, but also more difficult to analyze and solve.
Constrained nonlinear optimization problems on $\R^m$ can be cast in the following form,
\begin{equation}
	\label{eq:NLP}
	\begin{aligned}
		\text{Minimize}
		\quad
		&
		f(x)
		,
		\quad
		\text{where }
		x \in \R^m
		\\
		\text{subject to (\st)}
		\quad
		&
		g(x) \in K
		.
	\end{aligned}
\end{equation}
Here $f \colon \R^m \to \R$ denotes the objective function, and $g \colon \R^m \to \R^n$ represents the constraint function.
Moreover, $K \subset \R^n$ is a convex cone satisfying $0 \in K$, \ie, it induces a preorder on $\R^n$ defined by
\begin{equation*}
	y \le_K z
	\quad
	\Leftrightarrow
	\quad
	y - z \in K
	.
\end{equation*}
The constraint in \eqref{eq:NLP} can thus be written as $g(x) \le_K 0$.

Problems of the form \eqref{eq:NLP} include classical nonlinear programming problems with equality and inequality constraints.
These are described by $g(x) = (g_I(x),g_E(x))^\transp$ and $K = \R_-^k \times \{0\}^{n-k} \subset \R^k \times \R^{n-k}$, where $\R_-^k$ is the non-positive orthant in $\R^k$.

It is well known that---under appropriate constraint qualifications---local minimizers of \eqref{eq:NLP} admit Lagrange multipliers, \ie, there exists $\mu \in \R^n$ such that
\begin{equation}
	\label{eq:NLP_KKT}
	\begin{aligned}
		&
		f'(x) + \sum_{i=1}^n \mu_i \, g_i'(x)
		=
		0
		,
		\\
		&
		\mu \in K^\circ
		\coloneqq
		\setDef{s \in \R^n}{s^\transp v \le 0 \text{ for all } v \in K}
		,
		\\
		&
		\mu^\transp g(x) 
		= 
		0
	\end{aligned}
\end{equation}
holds.
In short, we can write $f'(x) + \mu^\transp \, g'(x) = 0$ with $\mu \in K^\circ$ and $\mu^\transp g(x) = 0$.
The set $K^\circ$ is called the polar cone of~$K$.

\Cref{eq:NLP_KKT} is known as generalized Karush--Kuhn--Tucker (KKT) conditions pertaining to problem~\eqref{eq:NLP}.
We refer the reader to, \eg, \cite[Ch.~9]{Luenberger:1969:1}, \cite{ZoweKurcyusz:1979:1}, \cite{Troeltzsch:1984:1}, \cite[Ch.~5]{Jahn:2007:1}, \cite[Ch.~6]{Troeltzsch:2010:1}, for results in this direction in finite and infinite-dimensional spaces.

In this paper, we generalize \eqref{eq:NLP} to constrained optimization problems on manifolds, replacing $\R^m$ and $\R^n$ by finite-dimensional, smooth manifolds~$\cM$ and $\cN$, respectively.
Theory for the case of equality and inequality constraints $g \colon \cM \to \R^n$ has been considered in \cite{YangZhangSong:2014:1,BergmannHerzog:2019:1} and some algorithmic approaches have been discussed in \cite{LiuBoumal:2019:2,ObaraOkunoTakeda:2020:1}.
Theory and an algorithm for equality constraints of the form $g(p)=q_*$ with $g \colon \cM \to \cN$ were presented in \cite{SchielaOrtiz:2021:1}.
Here we aim to incorporate equality and inequality constraints for manifold-valued constraint mappings $g \colon \cM \to \cN$.

Such an extension is not straightforward since there is no natural way to define a cone (nor a preorder) on the manifold~$\cN$ which would take the role of the condition $g(x) \in K$.
We propose here to overcome this difficulty by requiring the constraint function to have values in a \emph{submanifold with corners}~$\cK \subset \cN$, a mathematical object that corresponds to a convex cone locally in adequate charts.

We thus consider the following class of problems,
\begin{equation}
	\label{eq:problem_setting_introduction}
	\begin{aligned}
		\text{Minimize}
		\quad
		&
		f(p)
		,
		\quad
		\text{where }
		p \in \cM
		\\
		\text{\st}
		\quad
		&
		g(p) \in \cK
		,
	\end{aligned}
\end{equation}
which generalizes \eqref{eq:NLP}.
The description of the feasible set as $\cF \coloneqq \setDef{p \in \cM}{g(p) \in \cK}$ turns out to be convenient and relevant in a number of situations.
Moreover, it will be shown that this description is independent of possibly varying parametrizations of the given problem.

\begin{figure}
	\centering
	\includegraphics[width=0.3\textwidth]{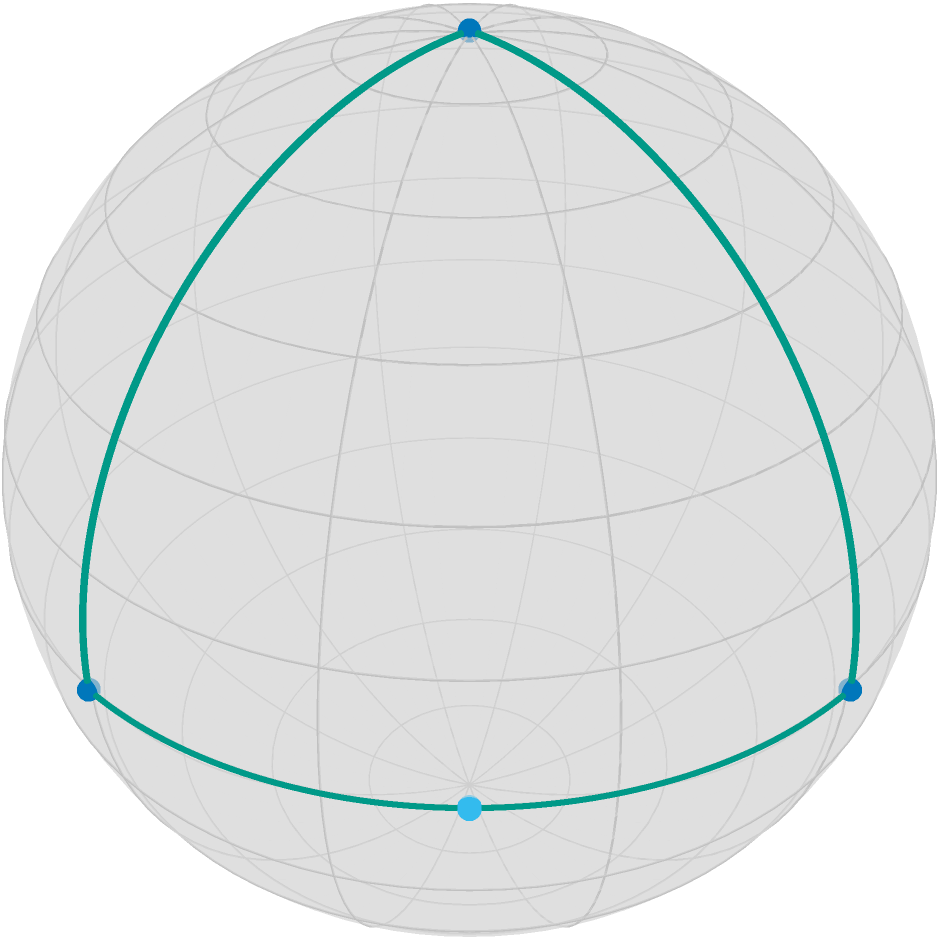}
	\caption{A geodesic polygon on the $2$-sphere. Unlike in $\R^2$, this set cannot be described as the intersection of half spaces. Notice that, for instance, at the tangent space at the light blue point in the middle of the horizontal geodesic, the image of the upper half space under the exponential map is the entire sphere.}
	\label{fig:geodesic_triangle}
\end{figure}

Our formulation differs from other generalizations of equality and inequality constraints.
Consider for instance a geodesic polygon as a feasible set~$\cF$, defined on the sphere $\cM = \sphere{2}$, \ie a set bounded by a set of geodesics.
More generally, we can also consider a geodesic polyhedron on $\sphere{m}$, \ie, a region bounded by a number of geodesic hyperplanes.
In other words, its boundary consists of totally geodesic submanifolds, \cf, \eg, \cite[Ch.~XI, §4]{Lang:1999:1}.
An example of a geodesic polygon is given in see \cref{fig:geodesic_triangle} in $\sphere{2}$.
$\cF$ constitutes a submanifold of $\cN = \cM$ with corners, so it can be naturally parametrized as $g(p) \in \cK$ with $g = \id_\cM$ and $\cK = \cF$.
By contrast, an algebraic description of $\cF$ in terms of classical inequalities runs into difficulties.
In the case of a vector space $\cM = \R^m$, the analogue of $\cF$ (an ordinary polygon) can be easily represented as the intersection of finitely many closed half spaces, using linear inequality constraints $g_i(x) = \inner{x - y_i}{n_i} \le 0$.
A similar attempt to describe $\cF$ on $\sphere{m}$ via inequality constraints of the type $g_i(p) = \riemannian{\logarithm{q_i} p}{n_i} \le 0$ can certainly be used locally; however, the lack of injectivity of the exponential map on $\sphere{m}$, and thus the lack of global well-definedness of its inverse, the logarithmic map, makes this inequality constraint globally not well-defined.

This paper is structured as follows.
We describe our approach to modeling manifold-valued constraints using manifolds with corners in \cref{section:constrained_optimization}.
Constraint qualifications are introduced and discussed in \cref{section:CQs}.
\Cref{section:first-order_conditions} is devoted to the derivation of first-order necessary optimality conditions.
We show in \cref{section:first-order_conditions_using_retractions} that equivalent conditions are obtained when the problem is pulled back to a tangent space, using a retraction.
In \cref{section:Lagrange_function} we introduce the analogue of a Lagrangian function for \eqref{eq:problem_setting_introduction}.
In preparation for the formulation of second-order optimality conditions in \cref{section:second-order_conditions}, we define the critical cone in \cref{section:Critical_cone}.
Finally, \cref{section:application_control_of_variational_problems} presents an application of our theory to the control of discretized variational problems.

We denote manifolds as well as subsets of manifolds by calligraphic letters.
For an introduction to differentiable manifolds, we refer the reader, \eg, to \cite{Lee:2012:1}.
Points on the manifold~$\cM$ are denoted by the letter~$p$, while points on $\cN$ are denoted by $q$.
Each manifold comes with a collection of charts $(\cU,\psi)$, and each chart maps an open subset $\cU$ of $\cM$ (or $\cN$) onto an open set in $\R^m$ (or $\R^n$), where $m$ and $n$ are the dimensions of $\cM$ and $\cN$, respectively.
We say that a chart $(\cU,\psi)$ is centered at a point $p$ if $p \in \cU$ holds.
For the purpose of this paper, since we will be pursuing a first- and second-order analysis, we will mostly assume that $\cM$ and $\cN$ are of class $C^2$, \ie, the chart transition maps $\psi_2 \circ \psi_1^{-1}$ are of this class.
In chart space, we use the letters $x \in \R^m$ and $y \in \R^n$.
We write $C^j(\cM,\cN)$ for the set of all mappings $\cM \to \cN$ which are $j$~times continuously differentiable.
The identity mappings on a vector space $V$ or on a manifold~$\cM$ are denoted by $\id_V$ and $\id_\cM$, respectively.
The zero element in the tangent space~$\tangentspace{p}{\cM}$ of a manifold~$\cM$ at $p$ is denoted by $0_p$.
We distinguish primal elements $v \in \tangentspace{p}{\cM}$ and dual elements $\mu \in \cotangentspace{p}{\cM}$ and write dual pairings in the form $\dual{\mu}{v}$ and compositions with linear mappings~$A$ into $\tangentspace{p}{\cM}$ as $\mu \, A$.

\section{Manifold-Valued Constraints}
\label{section:constrained_optimization}

Our method of choice to generalize equality and inequality constrained problems to manifolds is to replace the usual cone~$K$ that the equality and inequality constraints~$g$ are mapping into by a submanifold with corners.

In the following we use $0 \le k \le n$ and write $\R^k \times \{0\}^{n-k}$ to denote the subset of $\R^n$ consisting of those elements whose last $n-k$ components vanish.
We define the map $W \colon \R^n \to \R^{n-k}$ by $W x = (x_{k+1}, \dots, x_n)^\transp$.
Further, as usual, $v \le 0$ in $\R^\ell$ means $v_i \le 0$ for $i = 1, \dots, \ell$.

\begin{definition}[Submanifold with corners (\cite{Michor:1980:1})]
	\label{definition:submanifold_with_corners}
	Suppose that $\cN$ is an $n$-dimensional $C^2$-manifold.
	A subset $\cK \subset \cN$ is called a submanifold with corners of dimension~$k$ if, for each $q \in \cN$, there exists a local chart $(\cU,\psi)$ satisfying $\psi(q) = 0$, an index $\ell$ satisfying $0 \le \ell \le k$, and a \emph{surjective} linear operator
	\begin{equation*}
		A \colon \R^k \times \{0\}^{n-k} \to \R^\ell
	\end{equation*}
	such that
	\begin{equation*}
		\begin{aligned}
			\psi(\cK \cap \cU)
			&
			=
			\setDef{x \in \psi(\cU) \cap (\R^k \times \{0\}^{n-k})}{A \, x \le 0}
			\\
			&
			=
			\setDef{x \in \psi(\cU)}{A \, x \le 0, \; W x = 0}
		\end{aligned}
	\end{equation*}
	holds.
	In this case, $(\cU,\psi)$ is termed an adapted local chart centered at~$q$.
\end{definition}
We may identify $A$ with a matrix
$\begin{bsmallmatrix} \widehat A & 0 \end{bsmallmatrix}$
where $\widehat A \in \R^{\ell \times k}$ and $0 \in \R^{\ell \times (n-k)}$.
For $x \in \R^k \times \{0\}^{n-k}$, we then have $A \, x = \widehat A \, (x_1, \ldots, x_k)^\transp$.

We refer to $q$ in \cref{definition:submanifold_with_corners} as a corner of index~$\ell$.
It has been shown in \cite{Michor:1980:1} that the index~$\ell$, which may of course depend on $q$, however does not depend on the particular choice of the adapted local chart centered at~$q$.
In terms of optimization, $\ell$ describes the number of active inequality constraints at~$q$.
This generalizes the notion of vertices ($\ell = k$), edges ($\ell = k-1$), and higher-dimensional facets.

The requirement $\ell \le k$ is essential in this definition.
In local charts, the description of a corner satisfies the linear independence constraint qualification (LICQ), because the rows of $\widehat A$ are necessarily linearly independent to guarantee surjectivity.
Thus, whenever $(\widetilde \cU, \widetilde \psi)$ is a (non-adapted) local chart on $\cN$ such that $\widetilde \psi(\cK \cap \widetilde \cU)$ is given by the nonlinear constraint $\widetilde A(x) \le 0$ with $\widetilde A(0) = 0$, we can use the surjective implicit function theorem to construct an adapted local chart~$\psi$ such that $\psi(\cK \cap \cU)$ is described by $\widetilde A'(0) \, x \le 0$.

\Cref{definition:submanifold_with_corners} can be conceived as straightforward generalizations of the concepts
\begin{enumerate}
	\item
		of an embedded submanifold~$\cK \subset \cN$, which is obtained when $\ell = 0$ holds for all $q \in \cK$,
	\item
		of a smoothly bounded subset~$\cK \subset \cN$ with non-empty interior, which is obtained when $k = n$ and, for every $q \in \cN$, either $\ell = 0$ (interior point) or $\ell = 1$ (boundary point) holds,
	\item
		and of a convex polyhedron $\cK\subset \cN = \R^n$, whose corners satisfy the above regularity condition.
		In particular, the non-positive orthant~$\cK = \R_-^n\subset \R^n$ is a submanifold with corners of dimension~$n$ of $\R^n$.
		For instance, the origin $q = 0$ is a corner of index~$n$ and it can be described by $\widehat A = \id_{\R^n}$.
		As another example, the point $q = -e_j$ (the negative $j$-th unit vector in $\R^n$), is a corner of index~$n-1$ and a local description of $\cK$ can be defined via $\widehat A \in \R^{(n-1) \times n}$ whose rows are $e_i^\transp$ with $1 \le i \le n$, $i \neq j$.
\end{enumerate}

Next we discuss tangent spaces in the context of submanifolds with corners.
Among the various equivalent ways to define the tangent space for differentiable manifolds, we use the one given in \cite{Lang:1999:1,Michor:1980:1}.
Let $q \in \cN$ and consider the set
\begin{equation*}
	\setDef{(\psi,v)}{\psi \colon \cU_\psi \to \R^n \text{ is a chart at } q \in \cN, \; v \in \R^n}
	.
\end{equation*}
For two charts $\psi_1, \psi_2$, we denote the transition map by $T \coloneqq \psi_2 \circ \psi_1^{-1}$.
Define an equivalence relation $(\psi_1,v_{\psi_1}) \sim (\psi_2,v_{\psi_2})$ by
\begin{equation*}
	T'(\psi_1(q)) \, v_{\psi_1}
	=
	v_{\psi_2}
	.
\end{equation*}
We call any corresponding equivalence class a tangent vector $v$ of $\cN$ at~$q$ and $v_\psi$ its representative in the chart $\psi$.
For fixed $q \in \cN$, the set of these equivalence classes is a vector space $\tangentspace{q}{\cN}$, termed the tangent space of $\cN$ at~$q$.
The disjoint union of $\tangentspace{q}{\cN}$ over all $q \in \cN$ can be endowed with the structure of a manifold, more accurately a vector bundle, termed the tangent bundle $\tangentspace{}{\cN}$ of $\cN$.

Suppose now that $\cK$ is a submanifold with corners of $\cN$ of dimension~$k$.
For $q \in \cK$, we define the tangent space $\tangentspace{q}{\cK}$ as the set of all $v \in \tangentspace{q}{\cN}$ which possess a representative $v_\psi$ in an adapted chart $\psi$ centered at~$q$ such that $v_\psi$ is an element of $\R^k \times \{0\}^{n-k}$.
In this case, \emph{all} representatives of $v$ in all adapted charts centered at~$q$ satisfy the same relation.
It is easy to verify that $\tangentspace{q}{\cK}$ is a linear subspace of $\tangentspace{q}{\cN}$ of dimension~$k$.
Notice that the dimension of $\tangentspace{q}{\cK}$ does not depend on the index of $q$ as a corner of $\cK$.

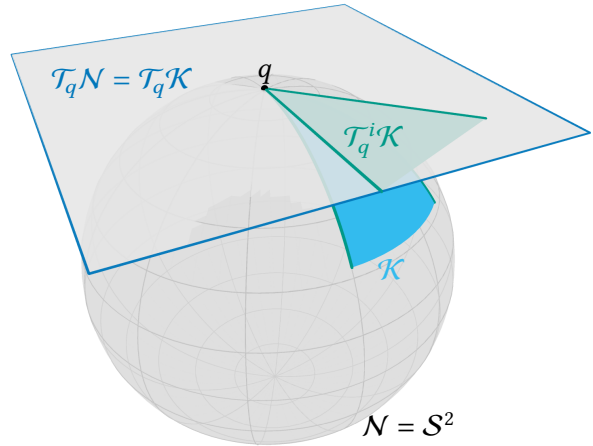
\begin{figure}
	\centering
	\begin{tikzpicture}
		\node[anchor = south west,inner sep = 0] (image) at (0,0) {\includegraphics[width = 0.5\textwidth]{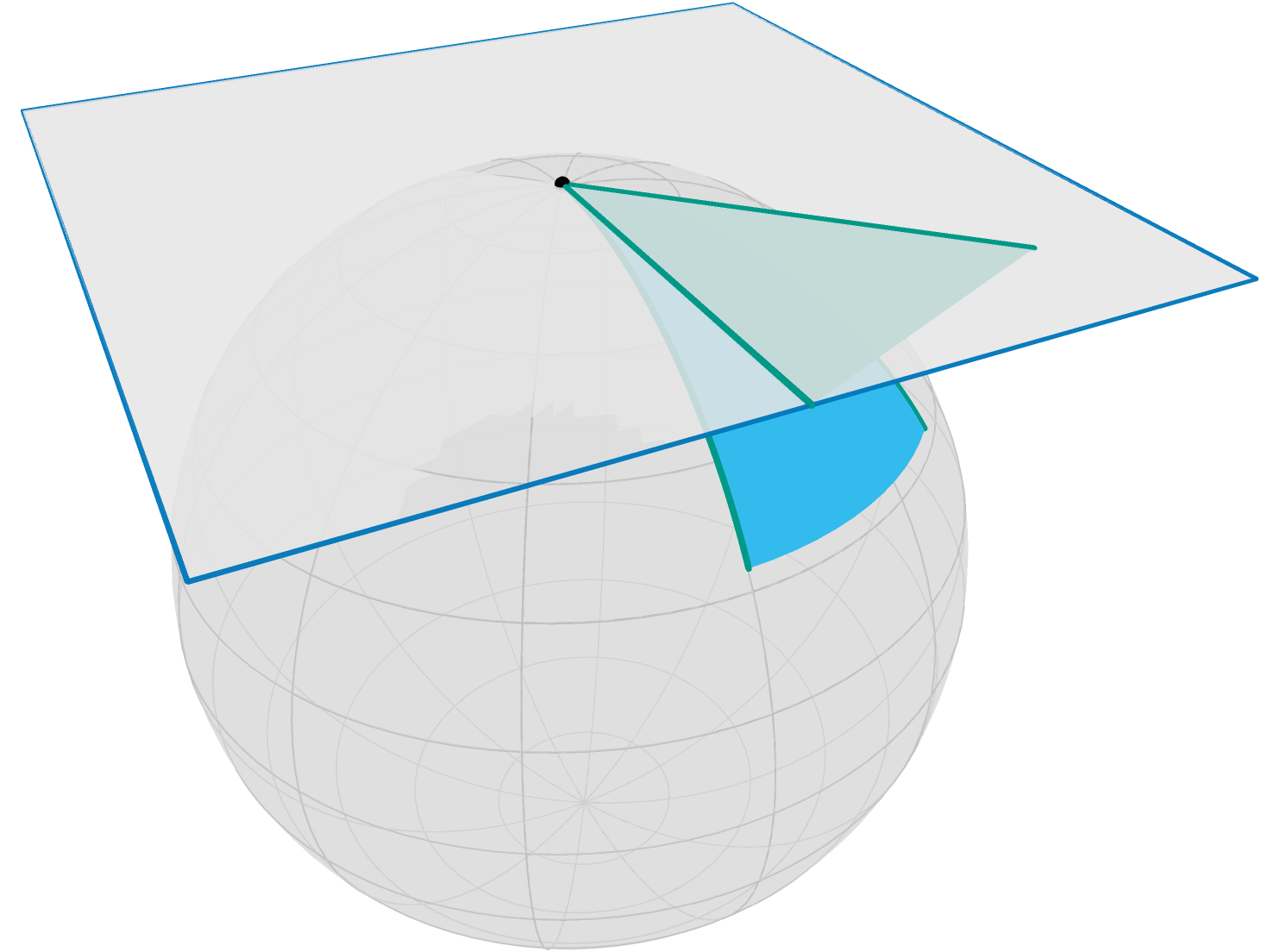}};
		\begin{scope}[x = {(image.south east)},y = {(image.north west)}]
			\node[text = TolVibrantBlue] at (0.2, 0.82) (TpM) {$\tangentspace{q}{\cN} = \tangentspace{q}{\cK}$};
			\node[text = TolVibrantTeal] at (0.62, 0.69) (TipK) {$\innertangentcone{q}{\cK}$};
			\node[text = TolVibrantCyan] at (0.65, 0.4) (KK) {$\cK$};
			\node at (0.68, 0.05) (NN) {$\cN = \sphere{2}$};
			\node at (0.44, 0.84) (q) {$q$};
		\end{scope}
	\end{tikzpicture}
	\caption{Illustration of a $k = 2$-dimensional manifold with corners~$\cK$ (teal) as a subset of the $n = 2$-dimensional sphere $\cN = \sphere{2}$. Due to $k = n$, the tangent space satisfies $\tangentspace{q}{\cK} = \tangentspace{q}{\cN}$ for every $q \in \cK$. At the particular point~$q$, which is a corner of index $\ell = 2$, the cone of inner tangent vectors $\innertangentcone{q}{\cK}$ is shown in green.}
	\label{figure:tangent_vectors}
\end{figure}

Further, the set of inner tangent vectors $\innertangentcone{q}{\cK} \subset \tangentspace{q}{\cK}$ is defined as all $v \in \tangentspace{q}{\cK}$ which satisfy, in addition, $A \, v_\psi \le 0$ for representatives in adapted charts centered at~$q$.
As discussed in \cite{Michor:1980:1}, $\innertangentcone{q}{\cK}$ is well-defined and it is a polyhedral convex cone.
Similarly, we denote by $\zerotangentspace{q}{\cK}$ the linear subspace of all elements~$v$ of $\tangentspace{q}{\cK}$ for which the representatives in adapted charts centered at~$q$ satisfy $A \, v_\psi = 0$.
We refer the reader to \cref{figure:tangent_vectors} for an illustrative example.

The following are our standing assumptions for the remainder of this paper.
\begin{assumption}
	Let $\cM$ and $\cN$ be $C^2$-manifolds of dimensions~$m$ and $n$, respectively.
	Moreover, let $\cK$ be a submanifold with corners of $\cN$ of dimension~$k$.
	We further suppose that $f \in C^2(\cM,\R)$ and $g \in C^2(\cM,\cN)$ hold and consider the following problem:
	\begin{equation}
		\label{eq:problem_setting}
		\begin{aligned}
			\text{Minimize}
			\quad
			&
			f(p)
			,
			\quad
			\text{where }
			p \in \cM
			\\
			\text{\st}
			\quad
			&
			g(p) \in \cK
			.
		\end{aligned}
	\end{equation}
\end{assumption}
Notice that products of submanifolds with corners are again submanifolds with corners.
One can therefore easily combine several constraints, \eg, $g_1(p) \in \cK_1$ and $g_2(p) \in \cK_2$, into one single constraint mapping into a product manifold.
We re-iterate that \eqref{eq:problem_setting} generalizes classical nonlinear programming problems with equality and inequality constraints.
The latter are obtained in case $\cM = \R^m$, $\cN = \R^n$, $\cK = \R_-^k \times \{0\}^{n-k} \subset \R^k \times \R^{n-k}$.
At any $p \in \cK$, the adapted local chart centered at a point $p$ can be chosen as $\varphi(\tilde p) = \tilde p-p$, and $\widehat A$ consists of the appropriate rows of $\id_{\R^k}$.

Be aware that in general the feasible set~$\cF \coloneqq g^{-1}(\cK) \subset \cM$ is \emph{not} a submanifold with corners even though $\cK$ is.
For example, consider $p$ to be the tip of a pyramid~$\cP$ in $\R^3$, where $\ell>3$ planes meet.
Then, locally near $p$, $\cP$ is described by $\ell>3$ inequality constraints, and thus $\cP$ cannot be a submanifold with corners of $\R^3$, because this would violate the condition $\ell \le k = 3$ in \cref{definition:submanifold_with_corners}.
Nevertheless, with a suitable affine mapping $g \colon \R^3 \to \R^\ell$, $\cP$ can be described locally as $\cP = g^{-1}(\R_-^\ell)$.
Thus, by means of the constraint mapping~$g$ we can obtain feasible sets more general than submanifolds with corners of $\cM$.
Also in view of practical computational approaches, the set~$\cK$ should have a simple structure, allowing, \eg, a local representation in computable adapted charts.

Suppose that $\varphi \colon \cM \supset \cU_p \to \R^m$ is a chart centered at~$p$ and that $\psi \colon \cN \supset \cU_{g(p)} \to \R^n$ is a chart centered at~$g(p)$.
We may then define the following local representations of $f$ and $g$:
\begin{equation*}
	f_\varphi
	\coloneqq
	f \circ \varphi^{-1}
	\colon
	\varphi(\cU_p)
	\to
	\R
	,
	\quad
	g_{\psi,\varphi}
	\coloneqq
	\psi \circ g \circ \varphi^{-1}
	\colon
	\varphi(\cU_p)
	\to
	\R^n
\end{equation*}
and obtain the following classical constrained optimization problem locally:
\begin{equation}
	\label{eq:problem_setting_local}
	\begin{aligned}
		\text{Minimize}
		\quad
		&
		f_\varphi(p_\varphi),
		\quad
		\text{where }
		p_\varphi \in \varphi(\cU_p)
		\\
		\text{\st}
		\quad
		&
		\paren[auto]\{.{%
			\begin{aligned}
				A \, g_{\psi,\varphi}(p_\varphi)
				&
				\le
				0
				,
				\\
				W g_{\psi,\varphi}(p_\varphi)
				&
				=
				0
				.
			\end{aligned}
		}
	\end{aligned}
\end{equation}
As a general strategy, we will carry over results on first- and second-order optimality conditions from \eqref{eq:problem_setting_local} to \eqref{eq:problem_setting} by formulations that are independent of the local representation in charts.
We will use rather straightforward and well established strategies of proof but highlight invariance considerations which arise in the differential geometric context.

\begin{example}\label{example:classical}
	Consider the standard case, \ie $\cN = \R^{n_I+n_E}$ and
	\begin{equation*}
		\begin{aligned}
			g_I(x)
			&
			\le
			0
			&
			&
			\text{in }
			\R^{n_I}
			,
			\\
			g_E(x)
			&
			=
			0
			&
			&
			\text{in }
			\R^{n_E}
			.
		\end{aligned}
	\end{equation*}
	This fits into our general setting \eqref{eq:problem_setting_introduction} if we define
	\begin{equation*}
		\begin{aligned}
			\cK
			\coloneqq
			\setDef[auto]{y \in \R^{n_I+n_E}}{%
				\begin{aligned}
					y_i
					&
					\le
					0
					\text{ for }
					i = 1, \ldots, n_I
					,
					\\
					y_i
					&
					=
					0
					\text{ for }
					i = n_I + 1,\ldots, n_I + n_E
				\end{aligned}
			}
			.
		\end{aligned}
	\end{equation*}
	This set $\cK$ is a submanifold of $\cN$ with corners of dimension~$k = n_I$.
	An adapted chart at a point $y \in \cK$ can be defined by $\varphi(\eta) = \eta - y$ and by choosing the chart domain~$\cU$ as an open $\norm{\cdot}_\infty$-ball about $y$ with radius $r = \min \setDef{\abs{y_i}}{y_i < 0}$.
	The index~$\ell$ of any point $y \in \cK$ equals the number of components~$1 \le i \le n_I$ for which $y_i = 0$ holds.
	Then the linear mapping $\widehat A \in \R^{\ell \times k}$ consists of rows equal to $e_i^\transp$ (the $i$-th unit vector in $\R^k$) for each index~$i$ with $y_i = 0$.
	For any $y \in \cK$, the tangent space $\tangentspace{y}{\cK}$ (in its representation \wrt the chart $\varphi(\eta) = \eta - y$) is given by $\R^{n_I} \times \{0\}^{n_E}$.
		At the point, $y = (1, 0, \ldots, 0)^\transp \in \cK$, for instance, the cone of inner tangent vectors is described by $\R \times \R_-^{n_I-1} \times \{0\}^{n_E}$, while the subspace $\zerotangentspace{y}{\cK}$ is equal to $\R \times \{0\}^{n_I-1} \times \{0\}^{n_E}$.
\end{example}

\begin{example}\label{eq:geopoly}
	Consider a geodesic polyhedron $\cK \subset \cN$ on a Riemannian manifold $\cN$, \ie, a set whose facets are totally geodesic submanifolds as in \cref{fig:geodesic_triangle}; \cf, \eg, \cite[Ch.~XI, §4]{Lang:1999:1}.
	We may use the logarithmic map $\logarithm{q} \colon \cN \to \tangentspace{q}{\cN} \cong \R^n$ to construct an adapted local chart at a point~$q \in \cK$.
	Then $\cK$ can be represented as $A \, v_\psi \le 0$ and $\cK$ is a manifold with corners, provided that $A$ (which depends on $q$ and $\psi$, of course) is surjective at any $q \in \cK$.
\end{example}

\begin{example}\label{eq:constraintlr}
	Given two mappings $g_\ell, g_r \colon \cM \to \cN$, consider the equality constraint
	\begin{equation*}
		g_\ell(p)
		=
		g_r(p)
		.
	\end{equation*}
	Since $\cN$ in general is not a vector space, this constraint cannot be written in the usual form $g_\ell(p) - g_r(p) = 0$.
	However, it can be formulated as $g(p) \in \cK$ via the mapping
	\begin{equation*}
		g
		\colon
		\cM
		\ni
		p
		\mapsto
		(g_\ell(p), g_r(p))
		\in
		\cN \times \cN
	\end{equation*}
	with $\cK = \setDef{(q_1,q_2) \in \cN \times \cN}{q_1 = q_2}$ the diagonal submanifold of $\cN \times \cN$.
\end{example}

\begin{example}
	Consider a vector bundle $\pi \colon \cN \to \cB$, where $\cB$ and $\cN$ are smooth manifolds and $\pi$ is a smooth surjective map.
	In fact, the total space~$\cN$ of a vector bundle is a manifold with special structure in the sense that, for each $q$ in the base manifold~$\cB$, the preimages $\pi^{-1}(q)$ (called fibres) are linear spaces; see, \eg, \cite[Ch.~III]{Lang:1999:1}.

	In applications, a constraint mapping $g \colon \cM \to \cN$ of the form
	\begin{equation*}
		g(p)
		=
		0_{\pi(g(p))}
	\end{equation*}
	arises frequently, in particular when $\cN = \tangentspace{}{\cB}$ or $\cN = \cotangentspace{}{\cB}$ is the tangent bundle or cotangent bundle over $\cB$, respectively.
	Since the mapping $q \mapsto 0_{\pi(q)}$ is well-defined and smooth on vector bundles, this constraint is of the form discussed in \cref{eq:constraintlr}.

	If the fibres $\pi^{-1}(q)$ of $\cN$ are equipped with preorder cones $K_q \subset \pi^{-1}(q)$, then also inequality constraints of the form
	\begin{equation*}
		g(p)
		\le
		0_{\pi(g(p))}
		,
		\quad
		\text{\ie,}
		\quad
		g(p) \in -K_{\pi(g(p))}
	\end{equation*}
	can be included under suitable assumptions on the choice of cones.
\end{example}

\section{Constraint Qualifications}
\label{section:CQs}

We recapitulate the definition of the tangent cone of a subset $\cF \subset \cM$ and generalize basic results, known for optimization problems on vector spaces, to the case of manifolds with corners.
We recall that $t_k \searrow 0$ denotes a sequence of strictly positive real numbers that converges to~$0$.
\begin{definition}[Tangent cone] \label{definition:tangent_cone}
	Let $p \in \cF$ and $(\cU,\varphi)$ be a chart centered at~$p$.
	A tangent vector~$v \in \tangentspace{p}{\cM}$ is said to belong to the tangent cone $\tangentcone{p}{\cF} \subset \tangentspace{p}{\cM}$ at~$p$ if there exists a representative $v_\varphi$ in the chart~$\varphi$ and sequences $t_k \searrow 0$ and $x_{\varphi,k} \in \R^m$ such that
	\begin{equation}\label{eq:tangential_sequence}
		x_{\varphi,k} \to v_\varphi
		\text{ and }
		t_k \, x_{\varphi,k} \in \varphi(\cF \cap \cU)
	\end{equation}
	holds.
	We then call $x_{\varphi,k}$ a feasible tangential sequence for $v_\varphi$.
\end{definition}

The following result shows that \cref{definition:tangent_cone} does not depend on the chosen chart.
Indeed, the tangent cone can alternatively be defined without the use of a chart; compare \cite[Def.~3.2]{BergmannHerzog:2019:1}.
\begin{lemma}\label{lemma:tangentcone_representative}
	Property~\eqref{eq:tangential_sequence} holds for one representative of $v \in \tangentspace{p}{\cM}$ if and only if it holds for every representative of $v$.
\end{lemma}
\begin{proof}
	Consider two local charts $\varphi_1$ and $\varphi_2$ centered at $p$ and their smooth transition map $T = \varphi_2 \circ \varphi_1^{-1}$, defined in a neighborhood~$\cU$ of $0 = \varphi_1(p) = \varphi_2(p) = T(0)$.
	Then the corresponding representatives $v_{\varphi_1}$ and $v_{\varphi_2}$ of a tangent vector $v \in \tangentspace{p}{\cM}$ are related by $v_{\varphi_2} = T'(0) \, v_{\varphi_1}$.
	By differentiability of $T$ we obtain (for sufficiently large $k$ so that $t_k \, x_{\varphi_1,k} \in \cU$):
	\begin{equation*}
		x_{\varphi_2,k}
		\coloneqq
		\frac{T(t_k \, x_{\varphi_1,k})-T(0)}{t_k} \to T'(0) \, v_{\varphi_1}
		=
		v_{\varphi_2}
	\end{equation*}
	for any pair of sequences $x_{\varphi_1,k} \to v_{\varphi_1}$ and $t_k \searrow 0$.
	Hence, $v_{\varphi_1}$ safisfies \eqref{eq:tangential_sequence} if and only if $v_{\varphi_2}$ does.
\end{proof}

Obviously, $\tangentcone{p}{\cF}$ is a cone and $0 \in \tangentcone{p}{\cF}$.
Furthermore, it is closed.
To see this, consider a sequence $v^i \in \tangentcone{p}{\cF}$ which converges to $v \in \tangentspace{p}{\cM}$ with $v \neq 0$.
Using a chart, we have sequences $t_k^i \searrow 0$ and $x_{\varphi,k}^i \to v_\varphi^i$.
From these, appropriate diagonal sequences can be chosen to verify $v \in \tangentcone{p}{\cF}$.

The following simple lemma can be proved as in the standard case:
\begin{lemma}\label{lemma:positive_derivative}
	Let $f \in C^1(\cM,\R)$ and assume that $p$ is a local minimizer of $f$ on a set~$\cF \subset \cM$.
	Then $f'(p) \, v \ge 0$ holds for all $v \in \tangentcone{p}{\cF}$.
\end{lemma}
\begin{proof}
	Consider $v \in \tangentcone{p}{\cF}$ and a corresponding tangential sequence $t_k \, v_k \in \cF$ with representatives $x_{\varphi,k}$.
	Then, by optimality, $t_k^{-1}(f_\varphi(t_k \, x_{\varphi,k}) - f_\varphi(0)) \ge 0$ holds for $k \in \N$ sufficiently large.
	Since $x_{\varphi,k} \to v_\varphi$ we obtain $f'(0) \, x_{\varphi,k} \to f'(0) \, v_\varphi$, but since $t_k^{-1}(f_\varphi(t_k \, x_{\varphi,k}) - f_\varphi(0)) - f'(0) \, x_{\varphi,k} \to 0$ by differentiability, this limit has to be non-negative.
\end{proof}

The following result shows that the tangent cone to a submanifold with corners has  a particularly simple structure since it agrees with the cone of inner tangent vectors defined in \cref{section:constrained_optimization}:
\begin{proposition}
	Suppose that $\cK$ is a submanifold with corners of $\cN$ and $q \in \cK$.
	Then
	\begin{equation*}
		\tangentcone{q}{\cK}
		=
		\innertangentcone{q}{\cK}
		.
	\end{equation*}
\end{proposition}
\begin{proof}
	Let $v \in \tangentspace{q}{\cK}$.
	Consider an adapted local chart $\psi$ of $\cK \subset \cN$, centered at~$q$, and defined on a neighborhood $\cU$ of $q$, and $v_\psi$ the corresponding representative of $v$.
	Since both $\tangentcone{q}{\cK}$ and $\innertangentcone{q}{\cK}$ are cones, we may assume \wolog that $\lambda \psi(\cK \cap \cU) \subset \psi(\cK \cap \cU)$ and $\lambda v_\psi \in \psi(\cU)$ for $\lambda \in [0,1]$.
	Two cases can occur.
	If $v_\psi \in \psi(\cK)$, then $v \in \innertangentcone{q}{\cK}$ holds by definition, and $v \in \tangentcone{q}{\cK}$ follows because $t_k \, v_\psi \in \psi(\cK \cap \cU)$ is clearly a tangential sequence.
	By contrast, if $v_\psi \not \in \psi(\cK)$, then $v \not \in \innertangentcone{q}{\cK}$ by definition.
	Moreover,
	\begin{equation*}
		\dist{\psi(\cK \cap \cU)}{v_\psi}
		\coloneqq
		\inf_{w \in \psi(\cK \cap \cU)} \norm{v_\psi - w}
		>
		0
	\end{equation*}
	because $\psi(\cK \cap \cU)$ is closed in $\psi(\cU)$.
	Then we can compute
	\begin{equation*}
		\dist{\psi(\cK \cap \cU)}{\lambda \, v_\psi}
		\coloneqq
		\inf_{w \in \psi(\cK \cap \cU)} \norm{\lambda(v_\psi - w)}
		=
		\lambda \dist{\psi(\cK \cap \cU)}{v_\psi}
		\;
		\text{for all } \lambda \in ]0,1]
		.
	\end{equation*}
	Hence, there is no feasible tangential sequence for $v_\psi$.
\end{proof}

In the following we consider the linearization
\begin{equation*}
	g'(p) \colon \tangentspace{p}{\cM} \to \tangentspace{g(p)}{\cN}
\end{equation*}
of $g$ at~$p$.
Its representation in a local chart $\varphi$, centered at $p$, and an adapted local chart $\psi$, centered at $g(p)$, reads:
\begin{equation*}
	g_{\psi,\varphi}'(0)
	\coloneqq
	(\psi \circ g \circ \varphi^{-1})'(0)
	\colon
	\R^m \to \R^n
	.
\end{equation*}

\begin{definition}[Linearizing cone]
	\label{definition:linearizing_cone}
	The linearizing cone at a point $p \in \cF$ is defined as
	\begin{equation*}
		\linearizingcone{p}{g}{\cK}
		\coloneqq
		\setDef[big]{v \in \tangentspace{p}{\cM}}{g'(p) \, v \in \innertangentcone{g(p)}{\cK}}
		=
		g'(p)^{-1}\paren[big](){\innertangentcone{g(p)}{\cK}}
		\subset
		\tangentspace{p}{\cM}
		.
	\end{equation*}
\end{definition}

\begin{lemma}
	We have $\tangentcone{p}{\cF} \subset \linearizingcone{p}{g}{\cK}$.
\end{lemma}
\begin{proof}
	Consider $v \in \tangentcone{p}{\cF}$, its representation in a chart $v_{\varphi}$ and corresponding sequences $t_k \searrow 0$ and $x_{\varphi,k} \to v_\varphi$, where $g_{\psi,\varphi}(t_k \, x_{\varphi,k}) \in \psi(\cK)$.
	We obtain:
	\begin{equation*}
		A \, g_{\psi,\varphi}(t_k \, x_{\varphi,k})
		\le
		0
		,
		\quad
		A \, g_{\psi,\varphi}(0)
		=
		0
		.
	\end{equation*}
	It follows that
	\begin{align*}
		t_k A \, g_{\psi,\varphi}'(0)(v_\varphi)
		&
		=
		A \, g_{\psi,\varphi}'(0)(t_k \, x_{\varphi,k}) + t_k A \, g_{\psi,\varphi}'(0)(v_\varphi - x_{\varphi,k})
		\\
		&
		=
		A \, g_{\psi,\varphi}(t_k \, x_{\varphi,k}) + \co(t_k)
		,
	\end{align*}
	and thus
	\begin{equation*}
		\frac{1}{t_k} A \, g_{\psi,\varphi}(t_k \, x_{\varphi,k}) \to A \, g_{\psi,\varphi}'(0)(v_\varphi)
		.
	\end{equation*}
	Since every row of the left hand side is non-positive, its limit cannot be positive.
	Thus, $A \, g_{\psi,\varphi}'(0)(v_\varphi) \le 0$ and similarly $W g_{\psi,\varphi}'(0)(v_\varphi) = 0$.
	This implies $v \in \linearizingcone{p}{g}{\cK}$.
\end{proof}

\begin{definition}
	The (description of the) feasible set $\cF$ is called transversal over $\cK$ at~$p \in \cF$ if
	\begin{equation*}
		\image g'(p) - \tangentspace{g(p)}{\cK}
		=
		\tangentspace{g(p)}{\cN}
		.
	\end{equation*}
	It is said to satisfy the Zowe--Kurcyusz--Robinson constraint qualification\\ (ZKRCQ, compare \cite{ZoweKurcyusz:1979:1}) at~$p \in \cF$ if
	\begin{equation}
		\label{eq:ZKRCQ}
		\tag{ZKRCQ}
		\image g'(p) - \innertangentcone{g(p)}{\cK}
		=
		\tangentspace{g(p)}{\cN}
		.
	\end{equation}
	It is said to satisfy the linear independence constraint qualification (LICQ) at~$p \in \cF$ if
	\begin{equation}
		\label{eq:LICQ}
		\tag{LICQ}
		\image g'(p) - \zerotangentspace{g(p)}{\cK}
		=
		\tangentspace{g(p)}{\cN}
		.
	\end{equation}
\end{definition}

Clearly, since $\zerotangentspace{g(p)}{\cK} \subset \innertangentcone{g(p)}{\cK} \subset \tangentspace{g(p)}{\cK}$ holds, \eqref{eq:LICQ} implies \eqref{eq:ZKRCQ}, which in turn implies transversality.
If the index $\ell$ of $g(p)$ satisfies $\ell = 0$, \ie $g(p)$ is not a corner of positive index, then all above notions are equivalent, because $\zerotangentspace{g(p)}{\cK} = \tangentspace{g(p)}{\cK}$ holds in this case.

\begin{proposition}\label{proposition:tangentcone_is_linearizingcone}
	If \eqref{eq:ZKRCQ} holds, then $\tangentcone{p}{\cF} = \linearizingcone{p}{g}{\cK}$.
\end{proposition}
\begin{proof}
	As above, consider a chart $\varphi$ of $\cM$ centered at~$p$ and an adapted chart $\psi$ of $\cN$ centered at~$g(p)$.
	Then the feasible set is represented locally as follows:
	\begin{equation*}
		\cF_{\psi,\varphi}
		\coloneqq
		\setDef{x \in \varphi(\cU)}{W g_{\psi,\varphi}(x) = 0, \; A \, g_{\psi,\varphi}(x) \le 0}
		,
	\end{equation*}
	while the representation of the linearizing cone is:
	\begin{equation}\label{eq:LRep}
		\linearizingcone{p}{g}{\cK}_{\psi,\varphi}
		\coloneqq
		\setDef{x \in \R^m}{W g_{\psi,\varphi}'(0) \, x = 0, \; A \, g_{\psi,\varphi}'(0) \, x \le 0}
		.
	\end{equation}
	Then \eqref{eq:ZKRCQ} can be written as:
	\begin{equation}\label{eq:ZKRCQ_Rep}
		\image g_{\psi,\varphi}'(0) - \setDef{y \in \R^n}{W y = 0, \; A y \le 0}
		=
		\R^n
		.
	\end{equation}
	Under assumption \eqref{eq:ZKRCQ}, we can apply \cite{ZoweKurcyusz:1979:1}
	to conclude that $\linearizingcone{p}{g}{\cK}_{\psi,\varphi}$ coincides with the tangent cone of $\cF_{\psi,\varphi}$ at~$0$ in $\R^n$, which is, by \cref{lemma:tangentcone_representative}, a representative of $\tangentcone{p}{\cF}$.
	Since both sets are representatives of subsets of $\tangentspace{p}{\cM}$, we conclude the result as claimed.
\end{proof}

Using the local representation \eqref{eq:LRep}, where our constraints are split into equality and inequality constraints, we can formulate the Mangasarian--Fromo\-vitz constraint qualification (MFCQ) in the following way:
\begin{equation}
	\label{eq:MFCQ}
	\tag{MFCQ}
	\paren[auto].\}{%
	\begin{aligned}
		&
		\text{The mapping $W g_{\psi,\varphi}'(0)$ is surjective.}
		\\
		&
		\text{There exists $\hat x \in \R^m$ such that $W g_{\psi,\varphi}'(0) \, \hat x = 0$}
		\\
		&
		\qquad
		\text{and $A \, g_{\psi,\varphi}'(0) \, \hat x < 0$ holds (in each component).}
	\end{aligned}
	\quad
	}
\end{equation}

\begin{proposition}
	\eqref{eq:MFCQ} and \eqref{eq:ZKRCQ} are equivalent.
\end{proposition}
\begin{proof}
	Let \eqref{eq:MFCQ} hold and $y \in \R^n$ be arbitrary.
	Define $\hat y \coloneqq g_{\psi,\varphi}'(0) \, \hat x \in \image g_{\psi,\varphi}'(0)$.
	In addition, since $W g_{\psi,\varphi}'(0)$ is surjective, there is $\tilde x$, such that $W g_{\psi,\varphi}'(0) \, \tilde x = W y$ and we define $\tilde y \coloneqq g_{\psi,\varphi}'(0) \, \tilde x$.
	Then we can write for any $\alpha > 0$:
	\begin{equation*}
		y
		=
		(\alpha \, \hat y + \tilde y) - (\alpha \, \hat y + \tilde y - y),
	\end{equation*}
	where $\alpha \, \hat y + \tilde y \in \image g_{\psi,\varphi}'(0)$.
	By construction, $W(\alpha \, \hat y + \tilde y - y) = 0$ holds, and choosing $\alpha$ sufficiently large we also obtain $A(\alpha \, \hat y + \tilde y - y) \le 0$, because $A \hat y < 0$.
	This shows \eqref{eq:ZKRCQ_Rep} and thus \eqref{eq:ZKRCQ}.

	If \eqref{eq:ZKRCQ} holds, then for any $y \in \R^n$ there is $\hat y \in \image g_{\psi,\varphi}'(0)$, such that $W \hat y = W y$ and $A \hat y \le Ay$, because $y = \hat y - (\hat y - y)$ with $W(\hat y - y) = 0$ and $A(\hat y - y) \le 0$.
	Thus, since $W$ and $A$ are surjective by definition of manifolds with corners, $W g_{\psi,\varphi}'(0)$ is surjective as well, and we find $y$ such that $W y = 0$ and $A y < 0$, and thus also $\hat y = g_{\psi,\varphi}'(0) \, \hat x$ with the same properties.
	So \eqref{eq:MFCQ} holds.
\end{proof}

\begin{proposition}\label{proposition:LICQ}
	$\cF$ satisfies \eqref{eq:LICQ} at $p \in \cF$ if and only if, for every representation in charts, the following linear mapping is surjective:
	\begin{equation*}
		B \, g_{\psi,\varphi}'(0)
		\colon
		\R^m \to \R^\ell \times \R^{n-k}
		,
		\quad
		\text{where }
		B
		\coloneqq
		\begin{pmatrix}
			A
			\\
			W
		\end{pmatrix}
		.
	\end{equation*}
\end{proposition}
\begin{proof}
	Let $v \in \tangentspace{g(p)}{\cN}$ with representative $v_\psi \in \R^n$.
	If $B \, g_{\psi,\varphi}'(0)$ is surjective, then we find $w_\varphi \in \R^m$, such that $B \, g_{\psi,\varphi}'(0) \, w_\varphi = -B \, v_\psi$.
	This implies that $v^0_\psi \coloneqq g_{\psi,\varphi}'(0) \, w_\varphi + v_\psi \in \ker B$ and we may write $v_\psi = g_{\psi,\varphi}'(0) \, w_\varphi - v^0_\psi$. Thus, we have found $w \in \tangentspace{p}{\cM}$ and $v^0 \in \zerotangentspace{g(p)}{\cK} $, such that $v = g'(p)w-v^0$.

	If, conversely, \eqref{eq:LICQ} holds, then we can write $v_\psi = g_{\psi,\varphi}'(0) \, w_\varphi - v^0_\psi$ for any $v \in \tangentspace{g(p)}{\cN}$ with $v^0_\psi \in \ker B$ and thus $B \, g_{\psi,\varphi}'(0) \, w_\varphi = B \, v_\psi$.
	Hence the surjectivity of $B \, g_{\psi,\varphi}'(0)$ follows from the surjectivity of $B$, which holds by \cref{definition:submanifold_with_corners} of a submanifold with corners.
\end{proof}

\section{First-Order Optimality Conditions}
\label{section:first-order_conditions}

In this section we address the first-order necessary optimality conditions for \eqref{eq:problem_setting} under the constraint qualification \eqref{eq:ZKRCQ}.
To this end, we recall that
\begin{equation*}
	S^\circ
	=
	\setDef{v^* \in V^*}{\dual{v^*}{s} \le 0 \text{ for all $s \in S$}}
\end{equation*}
denotes the polar cone of an arbitrary set~$S \subset V$ of a normed vector space~$V$.

\begin{theorem}
	\label{theorem:KKT}
	Suppose that $p_* \in \cF$ is a local minimizer of \eqref{eq:problem_setting} such that \eqref{eq:ZKRCQ} holds at~$p_*$.
	Then there exists a Lagrange multiplier $\mu \in \cotangentspace{g(p_*)}{\cN}$ such that the following KKT conditions hold:
	\begin{subequations}\label{eq:KKT_conditions}
		\begin{align}
			\label{eq:KKT_conditions1}
			&
			f'(p_*) + \mu \, g'(p_*)
			=
			0
			\quad
			\text{on }
			\cotangentspace{p_*}{\cM}
			,
			\\
			\label{eq:KKT_conditions2}
			&
			\mu
			\in
			\paren[big](){\innertangentcone{g(p_*)}{\cK}}^\circ
			.
		\end{align}
	\end{subequations}
	The set of all possible Lagrange multipliers,
		$\Lambda(p_*) = \setDef{\mu \in \cotangentspace{g(p_*)}{\cN}}{\eqref{eq:KKT_conditions} \text{ holds}}$
	is compact.
	If \eqref{eq:LICQ} holds, then $\Lambda(p_*)$ is a singleton.
\end{theorem}
\begin{proof}
	By \cref{lemma:positive_derivative} we have $f'(p_*) \ge 0$ on $\tangentcone{p_*}{\cF}$ and thus, by \cref{proposition:tangentcone_is_linearizingcone} on $\linearizingcone{p_*}{g}{\cK}$.
	Hence $v = 0$ is a minimizer of the following linear problem:
	\begin{equation*}
		\begin{aligned}
			\text{Minimize}
			\quad
			&
			f'(p_*) \, v
			,
			\quad
			\text{where }
			v \in \tangentspace{p_*}{\cM}
			\\
			\text{ \st}
			\quad
			&
			g'(p_*) \, v \in \innertangentcone{p_*}{\cK}
			.
		\end{aligned}
	\end{equation*}
	Due to the \eqref{eq:ZKRCQ} regularity condition, we can once more apply the results of \cite{ZoweKurcyusz:1979:1} to this problem to conclude the existence of a Lagrange multiplier $\mu$ such that the KKT conditions \eqref{eq:KKT_conditions} hold, so $\Lambda(p_*)$ is non-empty.
	Being the intersection of closed sets, $\Lambda(p_*)$ is also closed.

	In order to prove the boundedness of $\Lambda(p_*)$, we proceed by contradiction.
	Consider a sequence $\mu_k$ of Lagrange multipliers with $\norm{\mu_k} \to \infty$ and a corresponding bounded sequence $\lambda_k \coloneqq (\mu_k-\mu_1)/\norm{\mu_k}$ with $\mu_1/\norm{\mu_k} \to 0$.
	By picking a subsequence we may assume that $\lambda_k$ converges to a limit $\lambda_*$ with $\norm{\lambda_*} = 1$.
	Due to \eqref{eq:ZKRCQ}, every $v \in \tangentspace{g(p_*)}{\cN}$ can be written as $v = w-u$, where $w \in \image g'(p_*)$ and $u \in \innertangentcone{g(p_*)}{\cK}$.
	Then we compute
	\begin{equation*}
		\dual{\lambda_*}{v}
		=
		\lim_{k \to \infty} \dual{\lambda_k}{v}
		=
		\lim_{k \to \infty} \paren[auto](){\frac{1}{\norm{\mu_k}} \dual{(\mu_k-\mu_1)}{w} - \frac{\dual{\mu_k}{u}}{\norm{\mu_k}} + \frac{\dual{\mu_1}{u}}{\norm{\mu_k}}}
		.
	\end{equation*}
	Since $\dual{(\mu_k - \mu_1) }{ w} = 0$, $\dual{ \mu_k}{w} \le 0$, and the last addend in the sum tends to $0$, as $k \to \infty$, it follows that $\dual{ \lambda_*}{v} \ge 0$ holds for all $v \in \tangentspace{g(p_*)}{\cN}$ and thus $\lambda_* = 0$, which is in contradiction to $\norm{\lambda_*} = 1$.
	Hence, $\Lambda(p_*)$ is bounded and therefore compact.

	Now consider two solutions $\mu_1$ and $\mu_2$ of \eqref{eq:KKT_conditions}.
	Then $\mu_1 - \mu_2 \in \paren[big](){\zerotangentspace{g(p_*)}{\cK}}^\circ$ and $(\mu_1 - \mu_2) \, g'(p_*) = 0$.
	Hence, for all $v \in \image g'(p_*) - \zerotangentspace{g(p_*)}{\cK}$, it follows that $(\mu_1 - \mu_2) \, v = 0$.
	If \eqref{eq:LICQ} holds, then this implies $(\mu_1 - \mu_2) \, v = 0$ for all $v \in \tangentspace{g(p_*)}{\cN}$ and thus $\mu_1 = \mu_2$.
\end{proof}

In the following, we derive a representation $\mu_\psi \in \R^n$ of $\mu \subset \paren[big](){\innertangentcone{g(p)}{\cK}}^\circ$ with respect to an adapted local chart $\psi$ centered at~$g(p)$.
Recall that, by definition, $v \in \innertangentcone{g(p)} \cK$ holds if and only if $W v_\psi = 0$ and $A \, v_\psi \le 0$.

\begin{proposition}
	$\mu \in \paren[big](){\innertangentcone{g(p)} \cK}^\circ$ holds if and only if its representation $\mu_\psi$ in an adapted chart is of the following form:
	\begin{equation*}
		\mu_\psi
		=
		\begin{pmatrix}
			A^\transp \lambda_I
			\\
			W^\transp \lambda_E
		\end{pmatrix}
		\in \R^n
	\end{equation*}
	where $\lambda_I \ge 0 \in \R^\ell$ and $\lambda_E \in \R^{n-k}$.
	Hence, in local charts, \eqref{eq:KKT_conditions} reads:
	\begin{align*}
		f_\varphi'(p)
		+ \lambda_I^\transp A \, g_{\psi,\varphi}'(p)
		+ \lambda_E^\transp W g_{\psi,\varphi}'(p)
		&
		=
		0
		,
		\\
		\lambda_I
		&
		\ge
		0
		.
	\end{align*}
\end{proposition}
\begin{proof}
	Consider a representative $v_\psi$ of an element of $\innertangentcone{g(p)}{\cK}$ and $\mu_\psi$ of the claimed form:
	\begin{equation*}
		\dual{\mu_\psi}{v_\psi}
		=
		\dual{A^\transp \lambda_I}{v_\psi}
		+
		\dual{\lambda_E}{0}
		=
		\dual{\lambda_I}{A \, v_\psi}
		\le
		0
		.
	\end{equation*}
	Hence, $\mu \in \paren[big](){\innertangentcone{g(p)}{\cK}}^\circ$.

	For the converse, assume that $(\lambda_I)_i < 0$ for some $1\le i \le \ell$.
	Since $A$ is surjective, choose $v_\psi$ such that $A \, v_\psi = -e_i$ holds, which implies $\dual{A^\transp \lambda_I}{v_\psi} = -(\lambda_I)_i > 0$, so $\mu \not \in \paren[big](){\innertangentcone{g(p)}{\cK}}^\circ$.
\end{proof}

We return back to \cref{example:classical} and recall that the rows of $\widehat A \in \R^{\ell \times k}$ consist of those unit vectors $e_i^\transp$ for which $g_I(x_*)_i = 0$ holds.
We observe the representation
\begin{equation*}
	\mu_\psi
	=
	\begin{pmatrix}
		\eta_I
		\\
		\eta_E
	\end{pmatrix}
	,
	\quad
	\text{where }
	\eta_I
	=
	\sum_{\mrep{{\setDef{i}{g_I(x_*)_i = 0}}}{}} \lambda_i \, e_i
	\quad
	\text{with some }
	\lambda_i \ge 0
	.
\end{equation*}
Thus we obtain the classical complementarity result:
\begin{equation*}
	\eta_I \ge 0
	,
	\quad
	g_I(x_*) \le 0
	,
	\quad
	\dual{\eta_I}{g_I(x_*)}
	=
	0
	,
\end{equation*}
together with the well-known dual equation:
\begin{equation*}
	f'(x_*)
	+ \eta_I^\transp g_I'(x_*)
	+ \eta_E^\transp g_E'(x_*)
	=
	0
	.
\end{equation*}
After transposition, it takes the more familiar form
\begin{equation*}
	\nabla f(x_*)
	+ g_I'(x_*)^\transp \eta_I
	+ g_E'(x_*)^\transp \eta_E
	=
	0
	.
\end{equation*}

\section{Retractions and Linearizing Maps}
\label{section:first-order_conditions_using_retractions}

Numerical solution algorithms frequently employ retractions to pull back optimization problems on manifolds to the corresponding tangent spaces.
In this section we will consider reformulations of the KKT conditions \eqref{eq:KKT_conditions} in terms of these objects.
This is an alternative to our approach via local charts employed in \cref{section:first-order_conditions} and it allows us to argue more conveniently in some cases.
Moreover, retractions are also the approach we take for the second-order analysis in \cref{section:second-order_conditions}.

We will use the following definitions:
\begin{definition}\label{definition:retraction}
	Let $V_{0_p} \subset \tangentspace{p}{\cM}$ be a neighborhood of $0_p \in \tangentspace{p}{\cM}$.
	A $C^2$-mapping $\retract{p} \colon V_{0_p} \to \cM$ is called a local retraction at~$p$ if it satisfies:
	\begin{enumerate}[label=\ensuremath{(\retractionSymbol\roman*)},leftmargin=*]
		\item \label[property]{item:definition:retraction:zero}
			$\retract{p}(0_p) = p$,
		\item \label[property]{item:definition:retraction:diff}
			$D \retract{p}(0_p) = \id_{\tangentspace{p}{\cM}}$.
	\end{enumerate}

	Let $\cU_{q} \subset \cN$ be a neighborhood of $q \in \cN$.
	A $C^2$-mapping $S_q \colon \cU_q \to \tangentspace{q}{\cN}$ is called a local linearizing map at~$q$ if it satisfies:
	\begin{enumerate}[label=\ensuremath{(S\roman*)},leftmargin=*]
		\item \label[property]{item:definition:linearizing_map:y}
			$S_q(q) = 0_q$,
		\item \label[property]{item:definition:linearizing_map:Diff}
			$D S_q(q) = \id_{\tangentspace{q}{\cN}}$.
	\end{enumerate}
	We call $S_q$ adapted to $\cK$ if $S_q(\cU_q \cap \cK) = S_q(\cU_q) \cap \innertangentcone{q}{\cK}$ holds.
\end{definition}
Every chart $\varphi$ on $\cM$, centered at~$p$, induces a local retraction at~$p$ via $\retract{p}(v) \coloneqq \varphi^{-1}(v_\varphi)$.
Moreover, every adapted chart $\psi$ on $\cN$, centered at~$q$, induces an adapted linearizing map:
for any $\eta \in \cU_q$ we define $v \coloneqq S_q(\eta) \in \tangentspace{q}{\cN}$ by the equivalence class of $v_\psi \coloneqq \psi(\eta)$. If $\cK$ is a geodesic polyhedron on a Riemannian manifold $\cN$ as in \cref{eq:geopoly}, then $\logarithm{q}$ yields an adapted linearizing map at $q$.

\begin{remark}
	Retractions are widely used in optimization algorithms on manifolds; see, \eg, \cite{AbsilMahonySepulchre:2008:1}.
	Linearizing maps for constrained problems were introduced in \cite{SchielaOrtiz:2021:1}, but a similar concept has been used in a different context in~\cite{Boumal:2010:1} under the name \enquote{generalized logarithmic map}.

	The concept of adapted linearizing maps may be useful for the implementation of numerical algorithms in this setting.
	As we will see below, it allows us to write down a local optimization problem at~$p_*$ in a way that resembles a classical formulation without the need of further linearization of $S_{g(p_*)}(\cU_{g(p_*)} \cap \cK)$.
\end{remark}

Let $p$ be a feasible point of \eqref{eq:problem_setting} and $\retract{p} \colon V_{0_p} \to \cM$ and $S_{g(p)} \colon \cU_{g(p)} \to \tangentspace{g(p)}{\cN}$ be a given local retraction and adapted linearizing map, respectively.
Choosing their domain of definition sufficiently small, we may assume without loss of generality that $\retract{p}$ and $S_{g(p)}$ are injective with $g(\retract{p}(V_{0_p})) \subset \cU_{g(p)}$.
We can now locally pull back our problem as follows:
\begin{align*}
	\bf
	&
	\coloneqq
	f \circ \retract{p}
	\colon
	V_{0_p}
	\to \R
	,
	\\
	\bg
	&
	\coloneqq
	S_{g(p)} \circ g \circ \retract{p}
	\colon
	V_{0_p} \to \tangentspace{g(p)}{\cN}
	,
	\\
	\bK
	&
	\coloneqq
	S_{g(p)}(\cK \cap \cU_{g(p)}) = \innertangentcone{g(p)}{\cK} \cap S_{g(p)}(\cU_{g(p)})
	,
\end{align*}
and formulate a local optimization problem on the tangent space at~$p$:
\begin{equation}
	\label{eq:problem_setting_Retraction}
	\begin{aligned}
		\text{Minimize}
		\quad
		&
		\bf(v)
		,
		\quad
		\text{where }
		v \in V_{0_p} \subset \tangentspace{p}{\cM}
		\\
		\text{\st}
		\quad
		&
		\bg(v) \in \bK \subset \tangentspace{g(p)}{\cN}
		,
	\end{aligned}
\end{equation}
since $\bK$ is the intersection of a polyhedral convex cone and a neighborhood of $0_{g(p)}$.
It can thus be described by finitely many linear equality and inequality constraints on $\tangentspace{g(p)}{\cN}$.
Neglecting the local neighborhoods, \eqref{eq:problem_setting_Retraction} is locally a classical constrained optimization problem of the form:
\begin{equation}
	\label{eq:problem_setting_Retraction_classic}
	\begin{aligned}
		\text{Minimize}
		\quad
		&
		\bf(v)
		,
		\quad
		\text{where }
		v \in \tangentspace{p}{\cM}
		\\
		\text{\st}
		\quad
		&
		A_I \, \bg(v)
		\le
		0
		\\
		\text{and}
		\quad
		&
		A_E \, \bg(v)
		=
		0
	\end{aligned}
\end{equation}
with linear mappings $A_I \colon \tangentspace{g(p)} \cN \to \R^{\ell}$ and $A_E \colon \tangentspace{g(p)} \cN \to \R^{n-k}$.
{Notice that the data of problem \eqref{eq:problem_setting_Retraction_classic} is, of course, not uniquely defined.
For instance, we may premultiply $A_I$ by a positive diagonal matrix, and $A_E$ by any invertible matrix.
However, the viable choices for $A_I$ and $A_E$ do not depend on the choice of $S_{g(p)}$.}

\begin{theorem}\label{theorem:KKT_retractions}
	Suppose that $p_*$ is a feasible point of \eqref{eq:problem_setting}.
	Then $p_*$ is locally optimal for \eqref{eq:problem_setting} if and only if $v_* = 0 \in \tangentspace{p_*}{\cM}$ is a local minimizer of \eqref{eq:problem_setting_Retraction}.
	In this case, when \eqref{eq:ZKRCQ} holds at~$p_*$, then there exists $\mu \in \cotangentspace{g(p_*)}{\cN}$ such that
	\begin{align*}
		\bf'(0_{p_*}) + \mu \, \bg'(0_{p_*})
		&
		=
		0
		\quad
		\text{in }
		\cotangentspace{p_*}{\cM}
		,
		\\
		\mu
		\in
		\paren[big](){\tangentcone{0_{p_*}}{\bK}}^\circ
		&
		=
		\paren[big](){\innertangentcone{g(p_*)}{\cK}}^\circ
		.
	\end{align*}
\end{theorem}
\begin{proof}
	Clearly, $0_{p_*} \in \tangentspace{p_*}{\cM}$ is a local minimizer of
	\eqref{eq:problem_setting_Retraction} if and only if $p_*$ is a local minimizer of \eqref{eq:problem_setting}.
	Moreover, by the chain rule, using \cref{item:definition:retraction:diff} of $\retract{p_*}$ and \cref{item:definition:linearizing_map:Diff} of $S_{g(p_*)}$:
	\begin{equation*}
		\bf'(0_{p_*}) = f'(p_*)
		,
		\quad
		\bg'(0_{p_*}) = g'(p_*)
		,
		\quad
		\tangentcone{0_{p_*}}{\bK} = \innertangentcone{g(p_*)}{\cK}
		.
	\end{equation*}
	Thus, our conditions directly follow from \eqref{eq:KKT_conditions}.
\end{proof}

As an alternative approach, we can apply a classical theorem on KKT conditions to \eqref{eq:problem_setting_Retraction_classic} and obtain
\begin{equation}\label{eq:KKT_retractions_explicit}
	\begin{aligned}
		\bf'(0_{p_*}) + \lambda_I^\transp A_I \, \bg'(0_{p_*}) + \lambda_E^\transp A_E \, \bg'(0_{p_*})
		&
		=
		0
		,
		\\
		\lambda_I
		&
		\ge
		0
		,
	\end{aligned}
\end{equation}
with $\lambda_I \in \R^\ell$ and $\lambda_E \in \R^{n-k}$, which depend on the choice of $A_I$ and $A_E$.
By invariance, the first row equivalently yields:
\begin{equation*}
	f'(p_*) + \lambda_I^\transp A_I \, g'(p_*) + \lambda_E^\transp A_E \, g'(p_*)
	=
	0
\end{equation*}
and thus by comparison,
\begin{equation*}
	\mu
	=
	\lambda_I^\transp A_I + \lambda_E^\transp A_E
	\in
	\paren[big](){\innertangentcone{g(p_*)}{\cK}}^\circ
	.
\end{equation*}
We emphasize that the number of rows in $A_I$, which is equal to the index $\ell$ of the corner $g(p_*)$, depends on $g(p_*)$.
Thus, there is no further distinction necessary between active and inactive constraints, because this is already built into the local representation of $\cK$.

The formulation \eqref{eq:KKT_retractions_explicit} allows us to split the given constraints into individual components and to distinguish strongly active and weakly active constraints, according to the structure of $\lambda_I$.

\begin{definition}
	We call the $i$-th constraint $(A_I)_i \, \bg \le 0$ weakly active at $(p_*,\lambda_I,\lambda_E)$ if $(\lambda_I)_i = 0$ holds, and strongly active in case $(\lambda_I)_i > 0$.
\end{definition}
Observe that this definition does not depend on the particular choice of $A_I$.
If $A_I$ is premultiplied by a positive diagonal matrix, then the notion of weak and strong activity of $(A_I)_i$ is not changed.

\section{Lagrangian Functions}
\label{section:Lagrange_function}

When $\cN = V$ is a normed linear space with dual space $V^*$ and $g \colon \cM \to V$, then a Lagrangian function for our problem \eqref{eq:problem_setting} with Lagrange multiplier $\mu \in V^*$ can be defined as usual:
\begin{equation*}
	L
	\colon
	\cM \times V^*
	\ni
	(p,\mu)
	\mapsto
	L(p,\mu)
	\coloneqq
	f(p) + \mu(g(p))
	\in
	\R
	.
\end{equation*}
However when $\cN$ is a nonlinear manifold, then $\mu$ cannot be defined as a linear functional on $\cN$.
Rather, we need to replace it with a function $h \in C^1(\cN,\R)$ and define
\begin{equation*}
	L
	\colon
	\cM \times C^1(\cN,\R)
	\ni
	(p,h)
	\mapsto
	L(p,h)
	\coloneqq
	f(p) + h(g(p))
	\in \R
\end{equation*}
as a Lagrangian function.
In the following we will consider $h$ fixed and regard the mapping $p \mapsto L(p,h)\colon \cM \to \R$ as a function in $p$.
Its derivative $L'$ is given by
\begin{equation*}
	L'(p,h)
	\coloneqq
	\frac{\d}{\d p} L(p,h)
	=
	f'(p) + h'(g(p)) \, g'(p)
	.
\end{equation*}
For these derivatives to be well-defined at a point $p$, it is enough that $h$ is defined in some neighborhood of $p$.
We can observe two things.
First, $\mu \coloneqq h'(g(p)) \in \cotangentspace{g(p)}{\cN}$ can be interpreted as a Lagrange multiplier; second, $L'(p,h)$ only depends on $\mu = h'(g(p))$ and not on the particular choice of $h$.

The paragraph above explains how to obtain $\mu$ from $h$.
Conversely, let $p_* \in \cM$ be fixed and $q_* = g(p_*)$.
In view of the KKT-conditions \eqref{eq:KKT_conditions} we would like to extend a Lagrange multiplier $\mu\in \cotangentspace{q_*}{\cN}$ locally to a nonlinear function $h$ on a neighbourhood of $q_*$ such that $h'(q_*) = \mu$ holds.
This can be achieved by using a linearizing map $S_{q_*}$ about $q_*$ and defining $h \coloneqq \mu \circ S_{q_*}$.
Then we obtain a Lagrangian function of the form
\begin{equation*}
	L_{S_{q_*}}(p,\mu)
	\coloneqq
	L(p,\mu \circ S_{q_*})
	=
	f(p) + \mu \circ S_{q_*} \!\! \circ g(p)
	.
\end{equation*}
Since $h'(q_*) = \mu \circ DS_{q_*}(q_*) = \mu$, we obtain with this definition of $h$:
\begin{equation}
	\label{eq:coincidence_of_Lagrangian_functions}
	L'_{S_{q_*}}(p_*,\mu)
	=
	f'(p_*) + \mu \, g'(p_*) = L'(p_*,h)
	.
\end{equation}

Alternatively we may define Lagrangian functions near $p_*$ with $q_* = g(p_*)$ via pull-backs:
\begin{equation*}
	\begin{aligned}
		&
		\bL
		\colon
		\tangentspace{p_*}{\cM} \times \cotangentspace{q_*}{\cN}
		\to
		\R
		\\
		&
		(v,\mu)
		\mapsto
		\bL(v,\mu)
		\coloneqq
		\bf(v) + \mu(\bg(v))
		=
		(f \circ \retract{p_*})(v) + (\mu \circ S_{q_*} \!\! \circ g \circ \retract{p_*})(v)
	\end{aligned}
\end{equation*}
with derivative
\begin{equation*}
	\bL'(v,\mu)
	=
	\bf'(v) + \mu \, \bg'(v)
	\quad
	\text{and thus}
	\quad
	\bL'(0_{p_*},\mu)
	=
	f'(p_*) + \mu \, g'(p_*).
\end{equation*}
It is therefore justified to define the derivative of the Lagrangian function in the following way:
\begin{equation}
	\label{eq:coincidence_of_Lagrangian_functions_2}
	\begin{aligned}
		L'(p_*,\mu)
		\coloneqq
		f'(p_*) + \mu \, g'(p_*)
		=
		\bL'(0_{p_*},\mu)
		=
		L'_{S_{q_*}}(p_*,\mu)
		=
		L'(p_*,h)
		\\
		\text{for }
		\mu = h'(q_*)
		,
	\end{aligned}
\end{equation}
independently of the choice of the retraction $\retract{p_*}$, linearizing map~$S_{q_*}$, and $h$, as long as $\mu = h'(q_*)$.
Utilizing the identifications $\mu \coloneqq h'(g(p_*))$ and $h \coloneqq \mu \circ S_{g(p_*)}$, we find that the KKT conditions \eqref{eq:KKT_conditions} can equivalently be written in the familiar way:
\begin{subequations}\label{eq:KKT_conditions_pull-back}
	\begin{align}
		\label{eq:KKT_conditions_pull-back1}
		&
		L'(p_*,\mu)
		=
		0
		\quad
		\text{on }
		\cotangentspace{p_*}{\cM}
		,
		\\
		\label{eq:KKT_conditions_pull-back2}
		&
		\mu
		\in
		\paren[big](){\innertangentcone{g(p_*)}{\cK}}^\circ
		.
	\end{align}
\end{subequations}

\section{The Critical Cone}
\label{section:Critical_cone}

To derive second-order optimality conditions, we need a definition of the critical cone at a KKT point $p_*$ as a subset of the tangent cone~$\tangentcone{p_*}{\cF}$.
Suppose that $(p_*,\mu)$ satisfies the KKT conditions \eqref{eq:KKT_conditions}.
We define the critical cone at $p_*$ as
\begin{align*}
	\criticalcone{\cM}
	&
	\coloneqq
	\setDef{v \in \tangentspace{p_*}{\cM}}{g'(p_*) \, v \in \innertangentcone{g(p_*)}{\cK} \text{ and } f'(p_*) \, v = 0}
	\\
	&
	\mrep[r]{{}={}}{{}\coloneqq{}}
	\setDef{v \in \tangentspace{p_*}{\cM}}{g'(p_*) \, v \in \innertangentcone{g(p_*)}{\cK} \text{ and } \mu \, g'(p_*) \, v = 0}
	.
\end{align*}
We also introduce the definition
\begin{align*}
	\criticalcone{\cN}
	&
	\coloneqq
	g'(p_*) \, \criticalcone{\cM}
	=
	\setDef{w \in \innertangentcone{g(p_*)}{\cK}}{\dual{\mu}{w} = 0}
	\\
	&
	\mrep[r]{{}={}}{{}\coloneqq{}}
	\setDef{w \in \innertangentcone{g(p_*)}{\cK}}{(A_I)_j w = 0 \text{ for all } j = 1, \dots, \ell \text{ such that } (\lambda_I)_j = 0}
	,
\end{align*}
where $(A_I)_j$ are the components of the mapping $A_I \colon \tangentspace{g(p_*)}{\cN} \to \R^\ell$ used in \eqref{eq:problem_setting_Retraction_classic}.
Then we can write $\mu \in \paren[big](){\Span \criticalcone{\cN}}^\circ$ for any Lagrange multiplier $\mu \in \Lambda(p_*)$.

The following considerations will be useful for the discussion of second-order conditions:
\begin{lemma}\label{lemma:PhiCone}
	Suppose that $X$ is a normed linear space and $U$, $V$ are open neighborhoods of $0 \in X$.
	Consider a diffeomorphism $\Phi \colon U \to V$ such that $\Phi(0) = 0$ and $\Phi'(0) = \id_X$ hold.
	Let $K$ be a polyhedral cone of the form
	\begin{equation*}
		K
		=
		\setDef{v \in X}{A_I \, v \le 0, \; A_E \, v = 0}
	\end{equation*}
	with linear maps $A_I \colon X \to \R^{n_I}$ and $A_E \colon X \to \R^{n_E}$.
	Suppose that
	\begin{equation*}
		\Phi
		\colon
		K \cap U
		\to
		K \cap V
	\end{equation*}
	is bijective.
	Select a row $a_j = (A_I)_j$ and define the facet
	\begin{equation*}
		K_j
		=
		\setDef{v \in X}{A_I \, v \le 0, \; A_E \, v = 0, \; a_j v = 0}
		.
	\end{equation*}
	Then there are neighborhoods $\tilde U$ and $\tilde V$ of $0$ such that
	\begin{equation*}
		\Phi
		\colon
		K_j \cap \tilde U
		\to
		K_j \cap \tilde V
	\end{equation*}
	is also bijective.
\end{lemma}
\begin{proof}
	We may assume \wolog that $\tilde U = U = B_r(0)$ is an open ball of radius~$r$ about~$0$.
	Since $\Phi$ is a homeomorphism and thus preserves boundaries of sets, we conclude in particular that
	\begin{equation*}
		\Phi
		\colon
		\partial K \cap U
		\to
		\partial K \cap V
	\end{equation*}
	is also a homeomorphism.
	Consider now the \enquote{open} facet
	\begin{equation*}
		\tilde K_j
		=
		\setDef{v \in K_j}{(A_I)_\ell \, v < 0 \text{ for all } \ell \neq j}
		,
	\end{equation*}
	which is a relatively open subset of $\partial K$.
	Then $U \cap \tilde K_j$ is a connected set, because $U$ and $\tilde K_j$ are both connected and convex.
	The continuity of $\Phi$ implies that $\Phi(U \cap \tilde K_j)$ is connected as well.
	However, the arbitrary union of two (or more) distinct open facets is not connected because each $\tilde K_j$ is a relatively open subset of this union.
	Hence, $\Phi(U \cap \tilde K_j)$ is a subset of an open facet $\tilde K_\ell$ and it remains to show $j = \ell$.
	Since $\Phi'(0) = \id_X$ holds, we find that
	\begin{equation*}
		\Phi'(0)
		\colon
		\tilde K_j
		\to
		\tilde K_j
	\end{equation*}
	is bijective.
	Using the differentiability of~$\Phi$ this implies that there exists $x_0 \in \tilde K_j$ such that $\Phi(x_0) \in \tilde K_j$ holds.
	We thus conclude that $\Phi(U \cap \tilde K_j) \subset \tilde K_j$.

	Picking some $B_\rho(0) \subset V$ we can show by the same argumentation
	\begin{equation*}
		\Phi^{-1}(B_\rho(0) \cap \tilde K_j) \subset \tilde K_j \cap U
	\end{equation*}
	and thus $B_\rho(0) \cap \tilde K_j \subset \Phi(\tilde K_j \cap U)$.
	Thus, $\Phi(U \cap \tilde K_j)$ can be written as $\tilde K_j \cap \tilde V$, where $\tilde V$ is a neighborhood of $0$.
\end{proof}
This lemma can be applied recursively also to subfacets of $K$.
Hence, after finitely many steps of application, we conclude in particular that there are neighborhoods $U$ and $V$ of $0$ such that $\Phi$ maps $\criticalcone{\cN} \cap U$ bijectively onto $\criticalcone{\cN} \cap V$.

\begin{lemma}\label{lemma:ThetaIntoCone}
	Consider two adapted linearizing maps $S_{q,1}$ and $S_{q,2}$ and the transition map~$\Theta \coloneqq S_{q,1} \circ S_{q,2}^{-1}$.
	Then
	\begin{equation*}
		\begin{aligned}
			v \in \innertangentcone{q}{\cK}
			\quad
			&
			\Rightarrow
			\quad
			\Theta''(0_q)[v, v] \in \tangentspace{q}{\cK}
			,
			\\
			v \in \criticalcone{\cN}
			\quad
			&
			\Rightarrow
			\quad
			\Theta''(0_q)[v, v] \in \Span \criticalcone{\cN}
			.
		\end{aligned}
	\end{equation*}
\end{lemma}
\begin{proof}
	Consider any cone $K \subset \tangentspace{q}{\cN}$ such that $\Theta$ maps $K$ into $K$.
	Since $\Theta(0_q) = 0_q$ and $\Theta'(0_q) = \id_{\R^n}$ hold, we can compute
	\begin{equation*}
		\Theta''(0_q)[v, v]
		=
		\lim_{t \to 0} t^{-2} \paren[auto](){\Theta(t \, v) - \Theta(0_q) - \Theta'(0_q) \, t \, v}
		=
		\lim_{t \to 0} t^{-2} \paren[auto](){\Theta(t \, v) - t \, v}
		.
	\end{equation*}
	Since both $\Theta(t \, v)$ and $t \, v$ belong to $K$, $\Theta(t \, v) - t \, v$ belongs to $\Span K$ and thus so does the limit.
	By definition, $\Theta$ maps $\bK \subset \innertangentcone{q}{\cK}$ into $\innertangentcone{q}{\cK}$ and thus $\Theta''(0_q)[v,v] \in \Span \innertangentcone{q}{\cK} = \tangentspace{q}{\cK}$ for $v \in \innertangentcone{q}{\cK} $, proving our first assertion.
	Our second assertion follows similarly, because $\Theta$ maps $\criticalcone{\cN}$ into $\criticalcone{\cN}$ by \cref{lemma:PhiCone}.
\end{proof}

\section{Second-Order Optimality Conditions}
\label{section:second-order_conditions}

Compared to the case in vector spaces, the formulation of second-order conditions on manifolds exhibits an additional difficulty.
On a vector space $V$ the second derivative of a real-valued function $\sigma \colon V \to \R$ at $x \in V$ can be represented as a bilinear form $\sigma''(x) \colon V \times V \to \R$, whose definiteness properties can be studied.
In contrast, for $\sigma \colon \cM \to \R$ we have $\sigma' \colon \tangentspace{}{\cM} \to \R$ and thus $\sigma'' \colon \tangentspace{}{(\tangentspace{}{\cM})} \to \R$.
The required representation of $\sigma''(p)$ as a bilinear form on $\tangentspace{p}{\cM}$, \ie $\sigma''(p) \colon \tangentspace{p}{\cM} \times \tangentspace{p}{\cM} \to \R$, is not given canonically.
A connection, or equivalenty, a covariant derivative, has to be specified for this purpose.
However, at a stationary point $p_*\in \cM$, \ie $\sigma'(p_*) = 0$, second derivatives of scalar-valued functions \emph{can} be represented canonically by bilinear forms on $\tangentspace{p_*}{\cM}$ without the help of a covariant derivative, as shown in the following lemma.

\begin{lemma}\label{lemma:morse}
	Suppose that $\sigma \in C^2(\cM,\R)$.
	At a point~$p_* \in \cM$ satisfying $\sigma'(p_*) = 0$, the second derivative $\sigma''(p_*) \colon \tangentspace{p_*}{\cM} \times \tangentspace{p_*}{\cM} \to \R$ is a well-defined symmetric bilinear form, \ie, a symmetric $(2,0)$-tensor.
\end{lemma}
\begin{proof}
	Consider two charts ${\varphi_1}$ and ${\varphi_2}$ centered at $p_*$ so that $\varphi_1(p_*) = \varphi_2(p_*) = 0$ holds.
	Then $\sigma$ has representations $\sigma_{\varphi_1} \coloneqq \sigma \circ {\varphi_1}^{-1}$ and $\sigma_{\varphi_2} = \sigma \circ {\varphi_2}^{-1}$ in charts, and $\sigma_{\varphi_1} = \sigma_{\varphi_2} \circ T$ with $T = \varphi_2 \circ \varphi_1^{-1}$.
	Let $v_{\varphi_1}$ and $v_{\varphi_2}$ be the representatives of $v \in \tangentspace{p_*}{\cM}$.
	Then $v_{\varphi_2} = T'(0) \, v_{\varphi_1}$ holds and we have
	\begin{equation*}
		\sigma_{\varphi_1}'(0) \, v_{\varphi_1}
		=
		\sigma_{\varphi_2}'(0) \, T'(0) \, v_{\varphi_1}
		=
		\sigma_{\varphi_2}'(0) \, v_{\varphi_2}
		.
	\end{equation*}
	Using $\sigma'(p_*) = 0$ we find
	\begin{align*}
		\sigma_{\varphi_1}''(0)[v_{\varphi_1}, v_{\varphi_1}]
		&
		=
		\sigma_{\varphi_2}''(0)[T'(0) \, v_{\varphi_1}, T'(0) \, v_{\varphi_1}]
		+
		\sigma_{\varphi_2}'(0) \, T''(0)[v_{\varphi_1}, v_{\varphi_1}]
		\\
		&
		=
		\sigma_{\varphi_2}''(0)[T'(0) \, v_{\varphi_1}, T'(0) \, v_{\varphi_1}]
		.
	\end{align*}
	This implies the well-definedness of $\sigma''(p_*)$ on $\tangentspace{p_*}{\cM} \times \tangentspace{p_*}{\cM}$.
	Its symmetry follows from the theorem of Schwarz.
\end{proof}
As a consequence of \cref{lemma:morse}, second-order optimality conditions for unconstrained optimization problems on $C^2$-manifolds can be formulated without recourse to covariant derivatives.
Even for constrained problems for which the constraint target manifold $\cN = V$ is a \emph{linear space}, we can apply \cref{lemma:morse} to the Lagrangian function $L \colon \cM \times V^* \to \R$, \ie $\sigma(p) \coloneqq L(p,\mu)$, at a KKT point $p_*$ with Lagrange multiplier $\mu \in V^*$ and obtain a well-defined second derivative $L''(p_*,\mu) \colon \tangentspace{p_*}{\cM} \times \tangentspace{p_*}{\cM} \to \R$, because of $\sigma(p_*) = L'(p_*,\mu) = 0$.

For the general case of manifold-valued constraints, the situation is more complex, since, as we have seen, a classical Lagrange multiplier~$\mu$ cannot be used directly to define a Lagrangian function due to lack of linearity of $\cN$.
Instead, a nonlinear function $h \in C^2(\cN,\R)$ was used to define $L(p,h)$.
Although $L'(p,h)$ only depends on $\mu \coloneqq h'(g(p))$, the situation is different for the second-order derivative.
Let $p_*$ be a KKT-point, $q_* = g(p_*)$, and $\mu \in \cotangentspace{q_*}{\cN}$ the corresponding Lagrange multiplier such that $h'(q_*) = \mu$ and $L'(p_*,h) = 0$ hold.
Then we can apply \cref{lemma:morse} to $\sigma(p) \coloneqq L(p,h) = f(p) + h(g(p))$ and obtain a well-defined bilinear form at $p_*$:
\begin{equation*}
	L''(p_*,h)
	\colon
	\tangentspace{p_*}{\cM} \times \tangentspace{p_*}{\cM}
	\to
	\R
	.
\end{equation*}
Unfortunately, $L''(p_*,h)$ still depends on the particular choice of $h$ and not only on $\mu =h'(q_*)$.
This can be seen most clearly when $\cM$ and $\cN$ are linear spaces.
Then we can compute $L''(p_*,h)$ as follows:
\begin{equation*}
	L''(p_*,h)[v,v]
	=
	f''(p_*)[v,v]
	+ \mu \, g''(p_*)[v,v]
	+ h''(q_*)[g'(p_*) \, v, g'(p_*) \, v]
	,
\end{equation*}
and we observe that the third term on the right hand side depends on the second derivative of $h$. Of course, these second derivatives can be avoided when $\cN$ is a linear space by taking the canonical choice $h = \mu$, but such a canonical choice is not possible when $\cN$ is nonlinear.

However, suppose we use an \emph{adapted} linearizing map $S_{q_*}$ about $q_*$ to define $h = \mu \circ S_{q_*}$ and thus $L_{S_{q_*}}(p,\mu) = L(p,h)$ holds.
In that case, as we will show now, $L''(p_*, h)[v, v] = L_{S_{q_*}}''(p_*,\mu)[v,v]$ \emph{is} independent of the particular choice of $S_{q_*}$ on the critical cone, \ie, for $v \in \criticalcone{\cM}$.
This is all we need in order to formulate second-order optimality conditions in an invariant way.

\begin{proposition}\label{pro:invarianceLpp}
	Suppose that $p_*$ is a KKT point, $q_* = g(p_*)$ holds and $\mu \in \cotangentspace{q_*}{\cN}$ is a corresponding Lagrange multiplier so that \eqref{eq:KKT_conditions_pull-back} is satisfied.
	Let $S_{q_*,1}$ and $S_{q_*,2}$ be adapted linearizing maps about $q_*$.
	Then
	\begin{equation}\label{eq:adaptedLMInvariance}
		L_{S_{q_*,1}}''(p_*,\mu)[v, v]
		=
		L_{S_{q_*,2}}''(p_*,\mu)[v, v]
		\quad
		\text{for all }
		v \in \criticalcone{\cM}
		.
	\end{equation}
	In view of \eqref{eq:coincidence_of_Lagrangian_functions_2}, we therefore also refer to $L_{S_{q_*,i}}''(p_*,h_i)$ simply as $L''(p_*,\mu)$.
	Moreover, for any pullback with retraction $\retract{p_*}$ and adapted linearizing map $S_{q_*}$, the relation
	\begin{equation*}
		L''(p_*,\mu)[v, v]
		=
		\bL''(0_{p_*},\mu)[v, v]
		\quad
		\text{for all }
		v \in \criticalcone{\cM}
	\end{equation*}
	holds.
\end{proposition}
\begin{proof}
	Defining $\Theta \coloneqq S_{q_*,1} \circ S_{q_*,2}^{-1}$, we observe $\mu \circ S_{q_*,2} = \mu \circ \Theta \circ S_{q_*,1}$.
	Consequently, for $w \in \tangentspace{p_*}{\cM}$ and $p = \retract{p_*}(w)$, we have
	\begin{align*}
		L_{S_{q_*,2}}(p,\mu) - L_{S_{q_*,1}}(p,\mu)
		&
		=
		\mu \circ (\id_{\tangentspace{q_*}{\cN}} - \Theta) \circ S_{q_*,1} \circ g(p)
		\\
		&
		=
		\mu \circ (\id_{\tangentspace{q_*}{\cN}} - \Theta) \circ \bg(w)
		.
	\end{align*}
	The first derivatives read
	\begin{equation*}
		L_{S_{q_*,2}}'(p,\mu) \, v - L_{S_{q_*,1}}'(p,\mu) \, v
		=
		\mu \circ \paren[big](){\id_{\tangentspace{q_*}{\cN}} - \Theta'(\bg(w))} \circ \bg'(w) \, v
		.
	\end{equation*}
	Since $p_*$ is stationary, second derivatives of $L_{S_{q_*}}(p,\mu)$
	are well-defined and can be computed as follows, using the fact that $\Theta'(0_{q_*}) = \id_{\tangentspace{q_*}{\cN}}$ holds:
	\begin{align*}
		\MoveEqLeft
		\paren[big](){L_{S_{q_*,2}}''(p_*,\mu) - L_{S_{q_*,1}}''(p_*,\mu)}[v, v]\\
 		&
		=
		\mu \circ \paren[big](){\id_{\tangentspace{q_*}{\cN}} - \Theta'(0_{q_*})} \, \bg''(0_{p_*})[v, v]
		- \mu \circ \Theta''(0_{q_*})[\bg'(0_{p_*}) \, v, \; \bg'(0_{p_*}) \, v]
		\\
		&
		=
		- \mu \circ \Theta''(0_{q_*})[\bg'(0_{p_*}) \, v, \; \bg'(0_{p_*}) \, v]
		.
	\end{align*}
	For $v \in \criticalcone{\cM}$ we conclude $\bg'(0_{p_*}) \, v \in \criticalcone{\cN}$ and thus we find by \cref{lemma:ThetaIntoCone}, using that the linearizing maps are adapted:
	\begin{equation*}
		\Theta''(0_{q_*})[\bg'(0_{p_*}) \, v, \bg'(0_{p_*}) \, v]
		\in
		\Span \criticalcone{\cN}
		.
	\end{equation*}
	By stationarity and by definition of $\criticalcone{\cN}$, we infer $\restr{\mu}{\Span \criticalcone{\cN}} = 0$ and thus
	\begin{equation}\label{eq:coincideOnCcrit}
		\mu \circ \Theta''(0_{q_*})[\bg'(0_{p_*}) \, v, \bg'(0_{p_*}) \, v]
		=
		0
		\quad
		\text{for all }
		v \in \criticalcone{\cM}
		,
	\end{equation}
	which yields the desired result.
\end{proof}

\begin{remark}\label{rem:secondorderLM}
  The conclusion of \cref{pro:invarianceLpp} can be extended slightly beyond the class of adapted linearizing maps: let us call $S_{q_*,1}$ and $S_{q_*,2}$ second-order consistent, if their transition map $\Theta \coloneqq S_{q_*,1} \circ S_{q_*,2}^{-1}$ satisfies $\Theta''(0_{q_*}) = 0$.
	Clearly, \eqref{eq:coincideOnCcrit} holds for second-order consistent linearizing maps, even for all $v \in \tangentspace{p_*}{\cM}$.
	Hence, \eqref{eq:adaptedLMInvariance} extends to linearizing maps each of which is second-order consistent with some adapted linearizing map.
\end{remark}
\begin{remark}
	The restriction to \emph{adapted} linearizing maps in \cref{pro:invarianceLpp} is natural, taking into account the definition of a manifold with corners via adapted local charts.
	To illustrate that this restriction is also essential (up to \cref{rem:secondorderLM}), consider $\cM = \cN = \R^2$ with $p = (p_1,p_2)^\transp$, $f(p) = -p_1$, $g = \id_\cM$ and $\cK = \setDef{p \in \cM}{p_1 \le 0}$.
	Then $0$ is a local minimizer of $f$, $\criticalcone{\cM} = \criticalcone{\cN} = \setDef{v \in \cM}{v_1 = 0}$ hold, and
	\begin{equation*}
		0
		=
		L'(0,\mu) \, v
		=
		\begin{pmatrix} -v_1 \\ 0 \end{pmatrix} + \mu \, \begin{pmatrix} v_1 \\ v_2 \end{pmatrix}
		\quad
		\text{implies}
		\quad
		\mu
		=
		\begin{pmatrix} 1 \\ 0\end{pmatrix}
		.
	\end{equation*}
 Using the adapted linearizing map $S_{0,1} = \id_\cM$, we obtain $\mu \circ S_{0,1}(p) = p_1$ and $L_{S_{0,1}}''(0,\mu)[v, v] = 0$, but using the non-adapted linearizing map $S_{0,2}(p) \coloneqq (p_1 + \alpha \, p_2^2, p_2)$ would yield $\mu \circ S_{0,2}(p) = p_1 + \alpha \, p_2^2$ and $L_{0,2}''(0,\mu)[v, v] = 2 \, \alpha \, v_1^2$.

	In general, it is also not possible to extend \eqref{eq:adaptedLMInvariance} beyond $v \in \criticalcone{\cM}$.
	Using the adapted local linearizing map $S_{0,3}(p) \coloneqq (p_1 + p_1 \, p_2,p_2)^\transp$ (when $\abs{p_2} < 1$, then $p_1 + p_1 p_2 \ge 0 \Leftrightarrow p_1 \ge 0$), we obtain $\mu \circ S_{0,3}(p) = p_1 + p_1 \, p_2$ and thus $L_{S_{0,3}}''(0,\mu)[v, v] = 2 \, v_1 \, v_2$, which concides with $L_{S_{0,1}}''(0,\mu)[v, v] = 0$ on $\criticalcone{\cM}$ but not on all of $\R^2$.
\end{remark}

Having achieved an invariant definition of $L''$ on the critical cone, second-order optimality conditions for manifold-valued constraints can now be reduced to the classical vector-valued case.
Suppose that $p_* \in \cM$ is a KKT point with Lagrange multiplier $\mu$.
For any choice of retraction at~$p_*$ and adapted linearizing map at~$g(p_*)$, we consider the second derivative of the pullback $\bL''(0_{p_*},\mu)$, which---as we have seen---is invariant on the critical cone $\criticalcone{\cM}$.

Invoking well-known results from the literature, we obtain the following second-order sufficient optimality conditions:
\begin{theorem}
	Assume that $p_* \in \cM$ and $\mu \in \cotangentspace{g(p_*)}{\cN}$ satisfy the KKT conditions~\eqref{eq:KKT_conditions}.
	Moreover, suppose that
	\begin{equation*}
		L''(p_*,\mu)[v, v]
		>
		0
		\quad
		\text{holds for all }
		v \in \criticalcone{\cM} \setminus \{0_{p_*}\}
		.
	\end{equation*}
	Then $p_*$ is a strict local minimizer of problem~\eqref{eq:problem_setting}.
\end{theorem}
\begin{proof}
	It is clear that this result holds for $\bL''(0_{p_*},\mu)$ and thus, by invariance, it also holds for $L''(p_*,\mu)$; see, \eg, \cite[Thm.~12.6]{NocedalWright:2006:1}.
\end{proof}

Concerning second-order \emph{necessary} optimality conditions, a wide variety of constraint qualifications can be found in the literature (\cf, \eg, \cite{HaeserRamos:2019:1} and references therein), leading to second-order conditions of various strength.
We restrict our discussion here to the simplest case:
\begin{theorem}
	Assume that $p_* \in \cM$ is a local minimizer of problem~\eqref{eq:problem_setting} and that \eqref{eq:LICQ} holds at $p_*$.
	Then $p_*$ satisfies the KKT conditions~\eqref{eq:KKT_conditions} with some Lagrange multiplier $\mu \in \cotangentspace{g(p_*)}{\cN}$.
	Moreover,
	\begin{equation*}
		L''(p_*,\mu)[v, v]
		\ge
		0
		\quad
		\text{holds for all } v \in \criticalcone{\cM}
		.
	\end{equation*}
\end{theorem}
\begin{proof}
	It is clear that this result holds for $\bL''(0_{p_*},\mu)$ and thus, by invariance, it also holds for $L''(p_*,\mu)$; see, \eg, \cite[Thm.~12.5]{NocedalWright:2006:1}.
\end{proof}

\section{Application to the Control of Discretized Variational Problems}
\label{section:application_control_of_variational_problems}

Suppose that $\cY$ and $\cU$ are smooth manifolds and consider the following energy minimization problem, parametrized (or controlled) by $u$:
\begin{equation*}
	\text{Minimize}
	\quad
	E(y,u)
	,
	\quad
	\text{where }
	y \in \cY
	,
\end{equation*}
which we replace by its stationarity condition:
\begin{equation*}
	0^*_y
	=
	c(y,u)
	\coloneqq
	\partial_y E(y,u) \in \cotangentspace{y}{\cY}
	.
\end{equation*}
Such a situation occurs frequently in the infinite-dimensional context of variational problems, where occasionally $\cY$ and/or $\cU$ are nonlinear, smooth manifolds.
Also the principle of stationary action, which is applied, \eg, in classical mechanics, leads to problems of a similar form.
After discretization, a similar problem on finite-dimensional manifolds is obtained.

Using the control variable~$u$, an optimal control problem or a parameter identification problem may then be formulated as follows:
\begin{equation*}
	\begin{aligned}
		\text{Minimize}
		\quad
		&
		f(y,u)
		,
		\quad
		\text{where }
		(y,u) \in \cY \times \cU
		\\
		\text{\st}
		\quad
		&
		0^*_y
		=
		c(y,u)
		.
	\end{aligned}
\end{equation*}
A simple concrete example, which has been considered in \eg, \cite[Ch.~6]{OrtizLopez:2020:1}, is the optimal control of a static inextensible flexible rod. Here $y : [0,1]\to \R^3$ is the configuration of the rod, $u$ is an applied force, and $E(y,u)$ is the total energy of the rod. 
Inextensibility is modelled by requiring $y'(t) \in \mathbb S^2$ for all $t \in [0,1]$, the unit sphere in $\R^3$, which renders $\cY$ a nonlinear manifold. 
An appropriate objective function~$f$ may comprise the distance of $y$ to some desired configuration and a Tychonov term for $u$. 
For details we refer to \cite[Ch.~6]{OrtizLopez:2020:1} and \cite{SchielaOrtiz:2021:1}.

Setting $p \coloneqq(y,u)$, $\cM \coloneqq \cY \times \cU$, $\cN \coloneqq \cotangentBundle[\cY]$, and taking $\cK$ to be the zero-section of $\cotangentBundle[\cY]$, \ie, the pairs $(y,0^*_y)\in \cotangentBundle[\cY]$, which can be identified with $\cK = \cY$, we observe that this problem fits into our theoretical framework, where the constraint mapping is defined as follows:
\begin{equation*}
	g
	\colon
	\cY \times \cU \to \cotangentBundle[\cY]
	\ni
	p
	=
	(y,u)
	\mapsto
	g(p)
	\coloneqq
	(y,c(y,u))
	.
\end{equation*}
To formulate first-order optimality conditions, we calculate the derivative at a feasible point:
\begin{equation*}
	g'(p)
	=
	(\id_{\cM},c(y,u))'(y,u)
	\colon
	\tangentspace{y}{\cY} \times \tangentspace{u}{\cU}
	\to
	\tangentspace{(y,0^*_y)}{(\cotangentBundle[\cY])}
	.
\end{equation*}
At $0_y$ we can utilize the canonical splitting (a connection or covariant derivative is not required here) of the cotangent's tangent space
\begin{equation*}
	\tangentspace{(y,0^*_y)}{(\cotangentBundle[\cY])}
	\isomorphic
	\tangentspace{y}{\cY} \times \cotangentspace{y}{\cY}
\end{equation*}
into the tangent space of the base manifold and a fibre. This allows us to write $g'(p)$ as a pair:
\begin{align*}
	g'(p)
	\colon
	\tangentspace{y}{\cY} \times \tangentspace{u}{\cU}
	&
	\to \tangentspace{y}{\cY} \times \cotangentspace{y}{\cY}\\
	\delta p
	=
	(\delta y, \delta u)
	&
	\mapsto
	g'(p) \, \delta p
	=
	\paren[big](){\delta y, c'(y,u)(\delta y,\delta u)}
	\in
	\tangentspace{y}{\cY} \times \cotangentspace{y}{\cY}
\end{align*}
and the tangent space of $\cK$ as:
\begin{equation*}
	 \innertangentcone{(y,0^*_y)}{\cK}
	 =
	 \tangentspace{y}{\cY}
	 =
	 \tangentspace{y}{\cY} \times \{0^*_y\}
	 \subset
	 \tangentspace{y}{\cY} \times \cotangentspace{y}{\cY}.
\end{equation*}
Thus the linearized constraints can be split into two parts, the first of which is redundant:
\begin{equation*}
	g'(p) \, \delta p \in \innertangentcone{(y,0_y^*)}{\cK}
	\quad \Leftrightarrow \quad
	\delta y \in \tangentspace{y}{\cY}
	,
	\quad
	c'(y,u)(\delta y,\delta u)
	=
	0^*_y
	.
\end{equation*}
Constraint qualifications are fulfilled at $p$, provided that $\image g'(p) - \tangentspace{y}{\cY} \times \{0^*_y\} = \tangentspace{y}{\cY} \times \cotangentspace{y}{\cY}$ holds.
This is the case if and only if $c'(y,u) \colon \tangentspace{y}{\cY} \times \tangentspace{u}{\cU} \to \cotangentspace{y}{\cY}$ is surjective.

A Lagrange multiplier $\mu$ is an element of
\begin{equation*}
	(T^i_{(y,0^*_y)} \cK)^\circ
	=
	(\tangentspace{y}{\cY})^\circ
	=
	(\tangentspace{y}{\cY} \times \{0^*_y\})^\circ
	=
	\{0_y^*\} \times (\tangentspace{y}{\cY})^{**}
	=
	\{0_y^*\} \times \tangentspace{y}{\cY}
	,
\end{equation*}
where the last identity is the canonical identification of the bidual space with the primal space.
A Lagrange multiplier thus is a pair
\begin{equation*}
	\mu
	=
	(0_y^*, \lambda) \in \cotangentspace{y}{\cY} \times \tangentspace{y}{\cY}
	.
\end{equation*}
These splittings yield $\mu \, g'(p) \, \delta p = \paren[big](){0_y^*\delta y, \lambda c'(y,u)(\delta y,\delta u)} = \lambda c'(y,u)(\delta y,\delta u)$ and thus the KKT-conditions read
\begin{equation*}
	0
	=
	f'(y,u)(\delta y,\delta u) + \lambda \, c'(y,u)(\delta y,\delta u)
	\quad
	\text{for all }
	(\delta y,\delta u) \in \tangentspace{y}{\cY} \times \tangentspace{u}{\cU}
	.
\end{equation*}
Since $c(y,u) = \partial_y E(y,u)$ is a linear form on $\tangentspace{y}{\cY}$, $c'(y,u)$ can be interpreted as a bilinear form on $(\tangentspace{y}{\cY} \times \tangentspace{u}{\cU}) \times \tangentspace{y}{\cY}$ and we have (notice that $\partial_{yy} E(y,u)$ is well-defined by \cref{lemma:morse}, since $\partial_yE(y,u) = 0$ holds):
\begin{equation*}
	\lambda \, c'(y,u)(\delta y,\delta u)
	\!=\!
	(\partial_{y} E)'(y,u)(\lambda,\!\delta y,\!\delta u)
	\!=\!
	\partial_{yy} E(y,u)(\lambda,\!\delta y) +\partial_{yu} E(y,u)(\lambda,\!\delta u)
\end{equation*}
Then the KKT conditions read in more detail:
\begin{align*}
	\partial_y f(y,u) \, \delta y + \partial_{yy} E(y,u)(\lambda,\delta y)
	&
	=
	0
	\quad
	\text{for all } \delta y \in \tangentspace{y}{\cY}
	,
	\\
	\partial_u f(y,u) \, \delta u + \partial_{yu} E(y,u)(\lambda,\delta u)
	&
	=
	0
	\quad
	\text{for all } \delta u \in \tangentspace{u}{\cU}
	,
	\\
	\partial_y E(y,u) \, \delta y
	&
	=
	0
	\quad
	\text{for all } \delta y \in \tangentspace{y}{\cY}
	.
\end{align*}
To write down a Lagrangian function and second-order conditions, we need adapted linearizing maps on the zero section of $\cotangentspace{}{\cY}$ at a KKT-point $p_* = (y_*,u_*)$ with $q_* = g(p_*) = (y_*,0_{y_*}^*)$.
Utilizing the above splitting, these are those mappings $S_{q_*} \colon \cotangentspace{}{\cY} \to \tangentspace{y}{\cY} \times \cotangentspace{y}{\cY}$ which map the zero section $\cK = \cY$ to the first factor of the product, \ie $0_\eta \mapsto (\delta y(\eta),0_y)$.
For a specific example, consider a $C^2$-retraction $\retract{y_*} \colon T_{y_*} \cY \to \cY$ with derivative $D \retract{y_*}(v) \colon \tangentspace{{y_*}}{\cY} \to \tangentspace{\retract{y_*}(v)}{\cY}$.
Then an adapted linearizing map can be given as:
\begin{align*}
	S_{q_*} (y,w)
	\coloneqq
	(v,w \, D \retract{y_*}(v))
	,
	\quad
	\text{where }
	v
	=
	\retractionSymbol^{-1}_{y_*}(y) \in \tangentspace{{y_*}}{\cY}
	.
\end{align*}
Since $w \in \cotangentspace{y}{\cY}$ holds, it follows that $w \, D \retract{y_*}(v) \in \cotangentspace{y_*}{\cY}$, and $S_{q_*} (y,0^*_y) = (v,0^*_{y_*})$, as required.
With the help of this linearizing map, the Lagrange multiplier~$\mu$ can be extended locally to a function $h \in C^2(\cotangentspace{}{\cY},\R)$ as follows:
\begin{equation*}
	h (\eta,w)
	=
	\mu \, S_{q_*}(\eta, w)
	=
	0^*_{y_*} \, v + w \, D \retract{y_*}(v) \lambda
	=
	w \, D \retract{y_*}(v)
	\lambda
\end{equation*}
and thus the Lagrangian function near $p_*$ reads:
\begin{equation*}
	L_{S_{q_*}}(p,\mu)
	=
	f(y,u) + \partial_y E(y,u) \, D \retract{y_*}(v) \lambda
	,
	\quad
	v
	=
	\retractionSymbol^{-1}_{y_*}(y)
	.
\end{equation*}
Its first derivative at a feasible point, where $\partial_y E(y,u) = 0$ holds, is given by
\begin{equation*}
	L'_{S_{q_*}}(p,\mu)(\delta p)
	\!=\!
	f'(y,u)(\delta y,\delta u) + (\partial_y E)'(y,u) (D \retract{y_*}(v) \lambda,\delta y,\delta u),
	\;
	v
	=
	\retractionSymbol^{-1}_{y_*}(y)
	.
\end{equation*}
For a the KKT point $p_*$ we observe $L'_{S_{q_*}}(p_*,\mu) = 0$, since $D \retract{y_*}(0_{y_*}) = \id_{\tangentspace{y_*}{\cY}}$.

Since $\innertangentcone{(y,0_y^*)}{\cK}$ is a linear subspace in our setting, the critical cone $\criticalcone{\cM}$ is given as the preimage of $\criticalcone{\cN} = \tangentspace{y_*}{\cY} \times \{0^*_{y_*}\}$ under $g'(p_*)$, so it is the set
\begin{align*}
	\criticalcone{\cM}
	&
	=
	\setDef{(\delta y,\delta u)}{c'(y_*,u_*)(\delta y,\delta u) = 0}
	\\
	&
	=
	\setDef{(\delta y,\delta u)}{\partial_{yy} E(y_*,u_*)(v,\delta y) + \partial_{yu} E(y_*,u_*)(v,\delta u) = 0 \; \forall v \in \tangentspace{y_*}{\cY}}
	.
\end{align*}
Finally, the second derivative of the Lagrangian at $p_*$ is well-defined on $\criticalcone{\cM}$ and can, at least formally, be written as:
\begin{equation*}
	\begin{aligned}
		L_{S_{q_*}}''(y_*,u_*,\lambda)[\delta p, \delta p]
		=
		(f''(y_*,u_*)+ (\partial_yE)''(y_*,u_*)(\lambda))[(\delta y,\delta u),(\delta y,\delta u)]
		\\
		\text{for all }
		(\delta y,\delta u) \in \criticalcone{\cM}
		.
	\end{aligned}
\end{equation*}
As a consequence of the restriction $(\delta y,\delta u) \in \criticalcone{\cM}$ and the fact that $S_{q_*}$ is adapted, terms containing $DD \retract{p_*}$ are not present in this formula, which reflects \cref{pro:invarianceLpp}.

\section{Conclusion and Outlook}
\label{section:conclusion_outlook}

In this paper we have extended the analysis of optimization problems on manifolds from vector space-valued constraints to the much more flexible case of manifold-valued constraints.
We have seen that such problems arise naturally when constraints are formulated in a geometric way, and in the optimal control of variational problems on manifolds.
We generalized the polyhedric structure required for inequality constraints by using submanifolds with corners and adapted local charts.

First-order optimality conditions were derived, which directly generalize the known cases.
An appropriate definition of the Lagrangian function and the formulation of well-defined second-order optimality conditions, however, revealed the significance of the above mentioned polyhedric structure, reflected by the important role played by adapted linearizing maps.
We emphasize that in order to derive the theory, Riemannian metrics or connections were not needed.

Most of the stated results may be generalized to infinite-dimensional Banach manifolds.
However, we expect additional technical difficulties.
First, it seems to be an open problem how to generalize \cref{definition:submanifold_with_corners} to the infinite dimensional case, \ie, to define corners of infinite index $\ell = \infty$ in a useful way.
Second, already in infinite-dimensional Banach spaces, optimality conditions exhibit a couple of topologcal subtleties, which have to be tackled in the case of Banach manifolds, as well.

Further, algorithmic approaches for this class of optimization problems are still to be developed, even in the finite-dimensional setting.
An idea would be to extend SQP methods to this setting. 
At every iterate $x_k$ we perform a local pull-back of the given problem to tangent spaces, using retractions and adapted linearizing maps. 
Locally, we end up with a problem of the form~\eqref{eq:problem_setting_Retraction_classic}.
A QP step may then be computed for this pull-back, and an update can be defined via a retraction.
A detailed realization of this basic idea is, however, subject to future research.

\section*{Data Availability}

Data sharing not applicable to this article as no datasets were generated or analyzed during the current study.

\printbibliography

\end{document}